\documentclass{article}
\usepackage{graphicx,float,color,fancybox,shapepar,setspace,hyperref}
\usepackage{subfigure}
\usepackage{pgf,tikz}
\usetikzlibrary{arrows}
\newtheorem{Theorem}{Theorem}[section]
\newtheorem{fact}{Fact}[section]
\newtheorem{Lemma}{Lemma}[section]

\newtheorem{corollary}{Corollary}[section]

\def\square{\hbox{\vrule height8pt depth0pt
\vbox{\hrule width7.2pt\vskip7.2pt\hrule width7.2pt}\vrule
height8pt depth0pt}\smallskip}
\def\qed{\hfill\square}
\begin{document}
\baselineskip 0.18in

\title{Planar Ramsey Numbers of Four Cycles Versus Wheels}
\author{{Chen Yaojun}$^1$, {Miao Zhengke}$^2$ and {Zhou Guofei}$^1$\\
1. Department of Mathematics, Nanjing University, Nanjing 210093, China\\
2. Department of Mathematics, Jiangsu Normal University}

\date{}
\maketitle
\openup 1.5\jot
\begin{abstract}
For two given graphs $G$ and $H$ the planar Ramsey number $PR(G,H)$ is the smallest integer $n$ such that every planar graph $F$ on $n$ vertices either contains a copy of $G$, or its complement contains a copy of $H$. In this paper, we first characterize some structural properties of $C_4$-free planar graphs, and then we  determine all planar Ramsey numbers $PR(C_4, W_n)$, for $n\ge 3$.
\end{abstract}
\section{Introduction}
In this paper, all graphs are simple. Given two graphs $G$ and $H$, the Ramsey number $R(G, H )$ is the smallest integer $n$ such that every graph $F$ on $n$ vertices contains a copy of $G$, or its complement contains a copy of $H$. The determination of Ramsey numbers is an extremely difficult problem. In this paper, we are interested in planar Ramsey numbers. For two given graphs $G$ and $H$ the planar Ramsey number $PR(G,H)$ is the smallest integer $n$ such that every planar graph $F$ on $n$ vertices either contains a copy of $G$, or its complement contains a copy of $H$. The concept of planar Ramsey number was introduced by Walker \cite{walker} in 1969 and by Steinberg and Tovey \cite{MR1244935} in 1993, independently.

For planar Ramsey number, all pairs of complete graphs was determined in \cite{MR1244935}. Gorgol and Rucinski \cite{PR_Cycle} determined all pairs of cycles. 
By combining computer search with some theoretical results, A. Dudek and A. Rucinski \cite{dudek} compute most of the planar Ramsey numbers $PR(G_1,G_2)$, where each of $G_1$ and $G_2$ is a complete graph, a cycle or a complete graph without one edge.
 In \cite{zhouC3}, Zhou et. al.  determined $PR(C_3, W_n)$ for $n\ge 3$.  

In this paper, we first characterize some structural properties of $C_4$-free planar graphs, and then we  determine all planar Ramsey numbers $PR(C_4, W_n)$, for $n\ge 3$.

\section{Some concepts and notations}
Let $G=(V(G),E(G))$ be a graph and $G^c$ the complement of $G$. We define $\varepsilon(G)=|E(G)|$. 

Let $v$ be a vertex in $G$, the neighborhood of $v$, denoted by $N_G(v)$, is the vertex set  consisting of the vertices which are adjacent to $v$. We define $N_G[v]=N_G(v)\cup \{v\}$. We denote by $d_G(v)=|N_G(v)|$ the degree of $v$ in $G$. The maximum and minimum degrees in $G$, will be denoted by $\Delta(G)$ and $\delta(G)$ respectively.

Let $U\subseteq V(G)$, denote by $G[U]$  the subgraph induced by $U$ in $G$. The independence number, the connectivity and the minimum degree of $G$, are denoted by $\alpha (G)$, $k(G)$ and $\delta (G)$ respectively. 

Let $v$ be a vertex in $G$ and $H$ be a subgraph in $G$, we denote by $v+H$ the graph in which every vertex of $H$ is adjacent to $v$.
A {\it wheel} $W_n=\{x\}+C_n$ is a graph of order $n+1$, where $x$ is called the hub of the wheel, $C_{n}$ is a cycle of length $n$, and $x$ is adjacent to each vertices of $C_{n}$.   

A graph $G$ of order $n$ is said to be Hamiltonian if it contains an $n$-cycle; and $G$ is said to be pancyclic if $G$ contains cycles of length $k$, for all $k=3,4,\cdots,n$.

A planar graph which is embedded in a plane is called a plane graph, a face of length $k$ in $G$ is called a $k$-face, whose boundary has exactly $k$ edges, denote by $F(G), F_k(G)$ the set of faces  and the set of $k$-faces of $G$, respectively.  Let  $f_k$ be the number of $k$-faces in $G$.

Let $C$ be a cycle of a plane graph $G$, we call $C$ a separating cycle of $G$ if both the inside and outside of $C$ have at least one vertex.

Let $G$ be a plane graph, we denote by $\Gamma(G)$ the set of edges which are not covered by any triangle. The cardinality  of $\Gamma(G)$ is denoted by $\tau(G)$. 

Let $G$ be a plane graph. We can construct a new graph $G^*$ from $G$ as follows: the vertex set of $G^*$ consists of all the faces of lengths at least 5, and for each pair of vertices $f,g\in V(G^*)$, $f$ and $g$ are adjacent if and only if $f$ and $g$ have exactly one vertex or have exactly one edge in common. For convenience, we call $G^*$ the vertex-edge-dual of $G$.

If $G$ has $e_1$ vertices of degree $d_1$, $e_2$ vertices of degree $d_2$,$\cdots$, $e_k$ vertices of degree $d_k$, we will denote the degree sequence of $G$ by $d_1^{e_1}d_2^{e_2}\cdots d_k^{e_k}$, where $d_1\le d_2\le \cdots \le d_k$.

We define $\delta(n,C_4)=\max\{\delta(G)|G$ is a $C_4$ free planar graph $\}$; We  denote by $M(n,C_4)$ the maximum number of edges among all $C_4$-free planar graphs. In \cite{chen}, Zhou and Chen determined all the values of $M(n,C_4)$ for $n\ge 30$.

Our main result is the following two theorems:

\begin{Theorem}\label{thm20101011}
Let $A=\{30,36,39,42\}\cup\{k|k\ge 44\}$ and $B=\{k|10\le k\le 29\}\cup\{31, 32, 33, 34, 35, 37, 38, 40, 41, 43\}$ be two  integer sets. Then

\noindent{\rm (i)} $\delta(n,C_4)=4$, for each $n\in A$; 

\noindent{\rm (ii)} $\delta(n,C_4)=3$, for each $n\in B$;

\noindent{\rm (iii)} $\delta(n,C_4)=2$, for each $5\le n\le 9$.
\end{Theorem}

\begin{Theorem}\label{thm20111103} The planar Ramsey numbers of $C_4$ versus $W_n$ is as follows: 
\[
PR(C_4,W_n)=\left\{
\begin{array}{ll}
10, & {\rm if}\ n=3;\\
9, & {\rm if}\ n=6;\\
\\
n+4, & {\rm if}\ n\in\{k|7\le k\le 25\}\cup\{27,28,29,30,31,33,34,36,37,39\};\\
\\
n+5, & {\rm if}\ n\in \{4,5,26,32,35,38\}\cup\{k|k\ge 40\}.
\end{array}
\right.
\]
\end{Theorem}

\section[]{The min-max degrees in  $C_4$-free planar graphs}
The following result can be implied by the famous Euler's Formula on plane graphs, and we omit the proof here.
\begin{Theorem}\label{thm3}
Let $G$ be a $C_4$-free plane graph, then  $\varepsilon(G)=\frac{15}{7}(n-2)-\frac{2}{7}\tau(G)-\frac{3}{7}f_6-\frac{6}{7}f_7-\cdots -\frac{3(r-5)}{7}f_r\le \frac{15}{7}(n-2)$ (where $r$ is the maximum length of faces in $G$).
\end{Theorem}

\begin{corollary}\label{coro20101010}
If $G$ is a $C_4$-free planar graph of order $n$, then 

\noindent{\rm (i)} $\delta(G)\le 4$.

\noindent{\rm (ii)} $\delta(G)\le 3$, if $n\le 29$.
\end{corollary}

{\bf Proof.} (i) Let $G$ be a $C_4$-free planar graph of order $n$. Suppose on the contrary that $\delta(G)\ge 5$, which implies that the number of edges of $G$ is at least $\frac{5}{2}n>\frac{15}{7}(n-2)$, this contradicts Theorem \ref{thm3}.

(ii) If $\delta(G)\ge 4$, by (i), we have $\delta(G)= 4$, thus the number of edges is $2n\le \frac{15}{7}(n-2)$, this implies that $n\ge 30$, a contradiction.\qed

\begin{Lemma}\label{lem130420}\cite{chen}
 $M(n,C_4)=\lfloor15(n-2)/7\rfloor-\mu$
for $n \ge 30$, where  $\mu = 1$ if $n \equiv 3({\rm mod} 7)$ or $n = 32, 33, 37$, and $\mu= 0$ otherwise.
\end{Lemma}

By using the program PLANTRI by Brinkmann and McKay \cite{plantri}, we have checked the fact of the follwing three facts.

\begin{figure}[H]
\centering 
\subfigure[]{ \definecolor{ffffff}{rgb}{1,1,1}
\definecolor{qqqqff}{rgb}{0,0,1}
\begin{tikzpicture}[line cap=round,line join=round,>=triangle 45,x=0.45cm,y=0.45cm]
\clip(-2.92,-0.62) rectangle (5.02,6.18);
\draw (0.54,2.64)-- (1.46,2.78);
\draw (0.54,2.64)-- (1.66,1.74);
\draw (0.54,2.64)-- (0.52,1.56);
\draw (0.54,2.64)-- (-0.68,1.3);
\draw (0.54,2.64)-- (0.76,4.46);
\draw (1.46,2.78)-- (0.76,4.46);
\draw (1.46,2.78)-- (2.24,2.84);
\draw (1.46,2.78)-- (2.82,1.48);
\draw (1.46,2.78)-- (1.66,1.74);
\draw (1.66,1.74)-- (2.82,1.48);
\draw (1.66,1.74)-- (3.16,0.64);
\draw (1.66,1.74)-- (1.48,0.86);
\draw (1.66,1.74)-- (0.52,1.56);
\draw (0.52,1.56)-- (1.48,0.86);
\draw (0.52,1.56)-- (-1.22,0.14);
\draw (0.52,1.56)-- (-0.68,1.3);
\draw (-0.68,1.3)-- (-1.22,0.14);
\draw (-0.68,1.3)-- (-1.54,0.6);
\draw (-1.54,0.6)-- (0.76,4.46);
\draw (0.76,4.46)-- (0.82,5.98);
\draw (0.76,4.46)-- (2.24,2.84);
\draw (0.82,5.98)-- (2.24,2.84);
\draw (2.24,2.84)-- (4.8,-0.38);
\draw (2.24,2.84)-- (2.82,1.48);
\draw (2.82,1.48)-- (4.8,-0.38);
\draw (2.82,1.48)-- (3.16,0.64);
\draw (3.16,0.64)-- (4.8,-0.38);
\draw (3.16,0.64)-- (0.98,0.02);
\draw (3.16,0.64)-- (1.48,0.86);
\draw (1.48,0.86)-- (0.98,0.02);
\draw (-1.22,0.14)-- (1.48,0.86);
\draw (-1.22,0.14)-- (0.98,0.02);
\draw (-1.22,0.14)-- (-2.76,-0.4);
\draw (-1.22,0.14)-- (-1.54,0.6);
\draw (-1.54,0.6)-- (-2.76,-0.4);
\draw (-1.54,0.6)-- (0.82,5.98);
\draw (0.82,5.98)-- (-2.76,-0.4);
\draw (0.82,5.98)-- (4.8,-0.38);
\draw (4.8,-0.38)-- (-2.76,-0.4);
\draw (4.8,-0.38)-- (0.98,0.02);
\draw (0.98,0.02)-- (-2.76,-0.4);
\draw (-0.68,1.3)-- (0.76,4.46);
\begin{scriptsize}
\fill [color=qqqqff] (0.54,2.64) circle (1.5pt);
\fill [color=qqqqff] (1.46,2.78) circle (1.5pt);
\fill [color=ffffff] (1.66,1.74) circle (1.5pt);
\draw(1.66,1.74) circle (1.5pt);
\fill [color=qqqqff] (0.52,1.56) circle (1.5pt);
\fill [color=qqqqff] (-0.68,1.3) circle (1.5pt);
\fill[color=ffffff] (0.76,4.46) circle (1.5pt);
\draw(0.76,4.46) circle (1.5pt);
\fill [color=qqqqff] (2.24,2.84) circle (1.5pt);
\fill [color=qqqqff] (2.82,1.48) circle (1.5pt);
\fill [color=qqqqff] (3.16,0.64) circle (1.5pt);
\fill [color=qqqqff] (1.48,0.86) circle (1.5pt);
\fill [color=ffffff] (-1.22,0.14) circle (1.5pt);
\draw(-1.22,0.14) circle (1.5pt);
\fill [color=qqqqff] (-1.54,0.6) circle (1.5pt);
\fill [color=qqqqff] (0.82,5.98) circle (1.5pt);
\fill [color=ffffff] (4.8,-0.38) circle (1.5pt);
\draw(4.8,-0.38) circle (1.5pt);
\fill [color=qqqqff] (0.98,0.02) circle (1.5pt);
\fill [color=qqqqff] (-2.76,-0.4) circle (1.5pt);
\end{scriptsize}
\end{tikzpicture}

} 
\hspace{0.1in} 
\subfigure[]{ \definecolor{ffffff}{rgb}{1,1,1}
\definecolor{qqqqff}{rgb}{0,0,1}
\begin{tikzpicture}[line cap=round,line join=round,>=triangle 45,x=0.5cm,y=0.5cm]
\clip(-3.02,1.42) rectangle (4.38,7.74);
\draw (0.58,4.16)-- (1.4,3.54);
\draw (0.58,4.16)-- (0.58,2.98);
\draw (0.58,4.16)-- (-0.12,3.5);
\draw (0.58,4.16)-- (0.02,4.64);
\draw (0.54,5.74)-- (0.58,4.16);
\draw (0.54,5.74)-- (1.4,3.54);
\draw (1.4,3.54)-- (1.76,4.7);
\draw (1.4,3.54)-- (2.58,2.62);
\draw (1.4,3.54)-- (0.58,2.98);
\draw (0.58,2.98)-- (2.58,2.62);
\draw (0.58,2.98)-- (0.58,2.22);
\draw (0.58,2.98)-- (-1.44,2.58);
\draw (0.58,2.98)-- (-0.12,3.5);
\draw (-0.12,3.5)-- (-1.44,2.58);
\draw (-0.12,3.5)-- (-0.7,4.02);
\draw (-0.12,3.5)-- (0.02,4.64);
\draw (0.02,4.64)-- (-0.7,4.02);
\draw (0.02,4.64)-- (0,5.9);
\draw (0.02,4.64)-- (0.54,5.74);
\draw (0.54,5.74)-- (0,5.9);
\draw (0.54,5.74)-- (0.48,7.62);
\draw (0.54,5.74)-- (1.76,4.7);
\draw (1.76,4.7)-- (0.48,7.62);
\draw (1.76,4.7)-- (4.2,1.82);
\draw (1.76,4.7)-- (2.58,2.62);
\draw (2.58,2.62)-- (4.2,1.82);
\draw (2.58,2.62)-- (0.58,2.22);
\draw (0.58,2.22)-- (4.2,1.82);
\draw (0.58,2.22)-- (-2.9,1.82);
\draw (0.58,2.22)-- (-1.44,2.58);
\draw (-1.44,2.58)-- (-2.9,1.82);
\draw (-1.44,2.58)-- (-1.14,4.24);
\draw (-1.44,2.58)-- (-0.7,4.02);
\draw (-0.7,4.02)-- (-1.14,4.24);
\draw (-0.7,4.02)-- (0,5.9);
\draw (0,5.9)-- (-1.14,4.24);
\draw (0,5.9)-- (0.48,7.62);
\draw (0.48,7.62)-- (-1.14,4.24);
\draw (0.48,7.62)-- (-2.9,1.82);
\draw (0.48,7.62)-- (4.2,1.82);
\draw (4.2,1.82)-- (-2.9,1.82);
\draw (-2.9,1.82)-- (-1.14,4.24);
\begin{scriptsize}
\fill [color=qqqqff] (0.58,4.16) circle (1.5pt);
\fill [color=qqqqff] (1.4,3.54) circle (1.5pt);
\fill[color=ffffff] (0.58,2.98) circle (1.5pt);
\draw(0.58,2.98) circle (1.5pt);
\fill [color=qqqqff] (-0.12,3.5) circle (1.5pt);
\fill [color=qqqqff] (0.02,4.64) circle (1.5pt);
\fill [color=ffffff] (0.54,5.74) circle (1.5pt);
\draw(0.54,5.74) circle (1.5pt);
\fill [color=qqqqff] (1.76,4.7) circle (1.5pt);
\fill [color=qqqqff] (2.58,2.62) circle (1.5pt);
\fill [color=qqqqff] (0.58,2.22) circle (1.5pt);
\fill [color=ffffff] (-1.44,2.58) circle (1.5pt);
\draw(-1.44,2.58) circle (1.5pt);
\fill [color=qqqqff] (-0.7,4.02) circle (1.5pt);
\fill [color=qqqqff] (0,5.9) circle (1.5pt);
\fill [color=ffffff] (0.48,7.62) circle (1.5pt);
\draw(0.48,7.62) circle (1.5pt);
\fill [color=qqqqff] (4.2,1.82) circle (1.5pt);
\fill [color=qqqqff] (-2.9,1.82) circle (1.5pt);
\fill [color=qqqqff] (-1.14,4.24) circle (1.5pt);
\end{scriptsize}
\end{tikzpicture}

} 
\hspace{0.1in} 
\subfigure[]{
\definecolor{ffffff}{rgb}{1,1,1}
\definecolor{qqqqff}{rgb}{0,0,1}
\begin{tikzpicture}[line cap=round,line join=round,>=triangle 45,x=0.45cm,y=0.45cm]
\clip(-2.7,-1.38) rectangle (5.14,5.22);
\draw (1.14,1.62)-- (1.74,0.18);
\draw (1.14,1.62)-- (0.46,1.22);
\draw (1.14,1.62)-- (0.84,2.88);
\draw (1.14,1.62)-- (1.5,2.86);
\draw (1.14,1.62)-- (1.92,1.34);
\draw (1.74,0.18)-- (1.92,1.34);
\draw (1.74,0.18)-- (3.2,0.04);
\draw (1.74,0.18)-- (2.64,-0.5);
\draw (1.74,0.18)-- (-0.2,-0.42);
\draw (1.74,0.18)-- (-0.02,0.38);
\draw (1.74,0.18)-- (0.46,1.22);
\draw (0.46,1.22)-- (-0.02,0.38);
\draw (0.46,1.22)-- (0.06,1.72);
\draw (0.46,1.22)-- (0.84,2.88);
\draw (0.84,2.88)-- (0.06,1.72);
\draw (0.84,2.88)-- (1.24,5.02);
\draw (0.84,2.88)-- (1.5,2.86);
\draw (1.5,2.86)-- (1.24,5.02);
\draw (1.5,2.86)-- (2.44,1.86);
\draw (1.5,2.86)-- (1.92,1.34);
\draw (1.92,1.34)-- (2.44,1.86);
\draw (1.92,1.34)-- (3.2,0.04);
\draw (3.2,0.04)-- (2.44,1.86);
\draw (3.2,0.04)-- (4.84,-1.08);
\draw (3.2,0.04)-- (2.64,-0.5);
\draw (2.64,-0.5)-- (4.84,-1.08);
\draw (2.64,-0.5)-- (-2.34,-1.08);
\draw (2.64,-0.5)-- (-0.2,-0.42);
\draw (-0.2,-0.42)-- (-2.34,-1.08);
\draw (-0.2,-0.42)-- (-0.96,0.5);
\draw (-0.2,-0.42)-- (-0.02,0.38);
\draw (-0.02,0.38)-- (-0.96,0.5);
\draw (-0.02,0.38)-- (0.06,1.72);
\draw (0.06,1.72)-- (-0.96,0.5);
\draw (0.06,1.72)-- (1.24,5.02);
\draw (1.24,5.02)-- (-0.96,0.5);
\draw (1.24,5.02)-- (-2.34,-1.08);
\draw (1.24,5.02)-- (4.84,-1.08);
\draw (1.24,5.02)-- (2.44,1.86);
\draw (2.44,1.86)-- (4.84,-1.08);
\draw (4.84,-1.08)-- (-2.34,-1.08);
\draw (-2.34,-1.08)-- (-0.96,0.5);
\begin{scriptsize}
\fill [color=qqqqff] (1.14,1.62) circle (1.5pt);
\fill [color=ffffff] (1.74,0.18) circle (1.5pt);
\draw(1.74,0.18) circle (1.5pt);
\fill [color=qqqqff] (0.46,1.22) circle (1.5pt);
\fill [color=qqqqff] (0.84,2.88) circle (1.5pt);
\fill [color=qqqqff] (1.5,2.86) circle (1.5pt);
\fill [color=qqqqff] (1.92,1.34) circle (1.5pt);
\fill [color=qqqqff] (3.2,0.04) circle (1.5pt);
\fill [color=qqqqff] (2.64,-0.5) circle (1.5pt);
\fill [color=qqqqff] (-0.2,-0.42) circle (1.5pt);
\fill [color=qqqqff] (-0.02,0.38) circle (1.5pt);
\fill [color=qqqqff] (0.06,1.72) circle (1.5pt);
\fill[color=ffffff] (1.24,5.02) circle (1.5pt);
\draw(1.24,5.02) circle (1.5pt);
\fill [color=qqqqff] (2.44,1.86) circle (1.5pt);
\fill [color=qqqqff] (4.84,-1.08) circle (1.5pt);
\fill [color=qqqqff] (-2.34,-1.08) circle (1.5pt);
\fill [color=qqqqff] (-0.96,0.5) circle (1.5pt);
\end{scriptsize}
\end{tikzpicture}
} 

\caption{3 non-isomorphic triangulations of 16 vertices with $\delta=5$.\label{fig20111024}}
\end{figure}
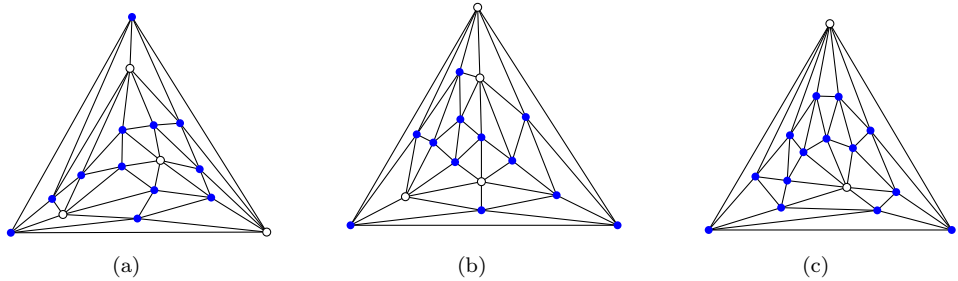

\begin{fact}\label{lem20111024}
There are all together 3 non-isomorphic triangulations of planar graphs on 16 vertices with minimum degree $\delta=5$ (Figure \ref{fig20111024}).
\end{fact}

\begin{figure}[H]
\centering 
\subfigure[]{\definecolor{ffffff}{rgb}{1,1,1}
\definecolor{qqqqff}{rgb}{0,0,0}
\begin{tikzpicture}[line cap=round,line join=round,>=triangle 45,x=0.3cm,y=0.3cm]
\clip(-4.42,-0.38) rectangle (3.63,6.86);
\draw (-0.22,2.71)-- (0.61,1.12);
\draw (-0.22,2.71)-- (-0.5,1.51);
\draw (-0.22,2.71)-- (-0.96,2.52);
\draw (-0.66,3.47)-- (-0.22,2.71);
\draw (-0.22,2.71)-- (0.12,3.22);
\draw (0.61,1.12)-- (0.12,3.22);
\draw (0.61,1.12)-- (1.55,1.3);
\draw (0.61,1.12)-- (1.11,0.38);
\draw (0.61,1.12)-- (-2.62,0.43);
\draw (0.61,1.12)-- (-0.5,1.51);
\draw (-0.5,1.51)-- (-2.62,0.43);
\draw (-0.5,1.51)-- (-1.47,1.56);
\draw (-0.5,1.51)-- (-0.96,2.52);
\draw (-0.96,2.52)-- (-1.47,1.56);
\draw (-0.96,2.52)-- (-1.74,2.41);
\draw (-0.96,2.52)-- (-0.66,3.47);
\draw (-0.66,3.47)-- (-1.74,2.41);
\draw (-0.66,3.47)-- (-0.45,5.01);
\draw (-0.66,3.47)-- (0.12,3.22);
\draw (0.12,3.22)-- (-0.45,5.01);
\draw (0.12,3.22)-- (0.42,3.82);
\draw (0.12,3.22)-- (1.55,1.3);
\draw (1.55,1.3)-- (0.42,3.82);
\draw (1.55,1.3)-- (3.4,-0.15);
\draw (1.55,1.3)-- (1.11,0.38);
\draw (1.11,0.38)-- (3.4,-0.15);
\draw (1.11,0.38)-- (-4.28,-0.15);
\draw (1.11,0.38)-- (-2.62,0.43);
\draw (-2.62,0.43)-- (-4.28,-0.15);
\draw (-2.62,0.43)-- (-2.69,1.93);
\draw (-2.62,0.43)-- (-1.47,1.56);
\draw (-1.47,1.56)-- (-2.69,1.93);
\draw (-1.47,1.56)-- (-1.74,2.41);
\draw (-1.74,2.41)-- (-2.69,1.93);
\draw (-0.45,5.01)-- (-1.74,2.41);
\draw (-0.45,5.01)-- (-2.69,1.93);
\draw (-0.45,5.01)-- (-0.48,6.74);
\draw (-0.45,5.01)-- (0.42,3.82);
\draw (0.42,3.82)-- (-0.48,6.74);
\draw (0.42,3.82)-- (3.4,-0.15);
\draw (3.4,-0.15)-- (-0.48,6.74);
\draw (3.4,-0.15)-- (-4.28,-0.15);
\draw (-4.28,-0.15)-- (-0.48,6.74);
\draw (-2.69,1.93)-- (-4.28,-0.15);
\draw (-2.69,1.93)-- (-0.48,6.74);
\begin{scriptsize}
\fill [color=qqqqff] (-0.22,2.71) circle (1.5pt);
\fill [color=ffffff] (0.61,1.12) circle (1.5pt);
\draw(0.61,1.12) circle (1.5pt);
\fill [color=qqqqff] (-0.5,1.51) circle (1.5pt);
\fill [color=qqqqff] (-0.96,2.52) circle (1.5pt);
\fill [color=qqqqff] (-0.66,3.47) circle (1.5pt);
\fill [color=ffffff] (0.12,3.22) circle (1.5pt);
\draw(0.12,3.22) circle (1.5pt);
\fill [color=qqqqff] (1.55,1.3) circle (1.5pt);
\fill [color=qqqqff] (1.11,0.38) circle (1.5pt);
\fill [color=ffffff] (-2.62,0.43) circle (1.5pt);
\draw(-2.62,0.43) circle (1.5pt);
\fill [color=qqqqff] (-1.47,1.56) circle (1.5pt);
\fill [color=qqqqff] (-1.74,2.41) circle (1.5pt);
\fill [color=ffffff] (-0.45,5.01) circle (1.5pt);
\draw (-0.45,5.01) circle (1.5pt);
\fill [color=qqqqff] (0.42,3.82) circle (1.5pt);
\fill [color=qqqqff] (3.4,-0.15) circle (1.5pt);
\fill [color=qqqqff] (-4.28,-0.15) circle (1.5pt);
\fill [color=ffffff] (-2.69,1.93) circle (1.5pt);
\draw (-2.69,1.93) circle (1.5pt);
\fill [color=qqqqff] (-0.48,6.74) circle (1.5pt);
\end{scriptsize}
\end{tikzpicture} }
\hskip 0.2in
\subfigure[]{ \definecolor{ffffff}{rgb}{1,1,1}
\begin{tikzpicture}[line cap=round,line join=round,>=triangle 45,x=0.31cm,y=0.33cm]
\clip(-2.76,0.22) rectangle (5.59,6.83);
\draw (1.44,2.48)-- (2,1.68);
\draw (1.44,2.48)-- (0.73,1.68);
\draw (1.44,2.48)-- (0.26,2.48);
\draw (1.44,2.48)-- (0.89,3.4);
\draw (1.44,2.48)-- (2.22,2.48);
\draw (2,1.68)-- (2.22,2.48);
\draw (2,1.68)-- (3.37,1.38);
\draw (2,1.68)-- (1.42,0.93);
\draw (2,1.68)-- (0.73,1.68);
\draw (0.73,1.68)-- (1.42,0.93);
\draw (0.73,1.68)-- (-0.46,1.4);
\draw (0.73,1.68)-- (0.26,2.48);
\draw (0.26,2.48)-- (-0.46,1.4);
\draw (0.26,2.48)-- (-0.58,2.48);
\draw (0.26,2.48)-- (0.85,4.75);
\draw (0.26,2.48)-- (0.89,3.4);
\draw (0.89,3.4)-- (0.85,4.75);
\draw (0.89,3.4)-- (2.08,3.4);
\draw (0.89,3.4)-- (2.22,2.48);
\draw (2.22,2.48)-- (2.08,3.4);
\draw (2.22,2.48)-- (3.14,2.46);
\draw (2.22,2.48)-- (3.37,1.38);
\draw (3.37,1.38)-- (3.14,2.46);
\draw (3.37,1.38)-- (5.41,0.42);
\draw (3.37,1.38)-- (1.42,0.93);
\draw (1.42,0.93)-- (5.41,0.42);
\draw (1.42,0.93)-- (-2.6,0.42);
\draw (1.42,0.93)-- (-0.46,1.4);
\draw (-0.46,1.4)-- (-2.6,0.42);
\draw (-0.46,1.4)-- (-0.58,2.48);
\draw (-0.58,2.48)-- (-2.6,0.42);
\draw (-0.58,2.48)-- (1.38,6.71);
\draw (-0.58,2.48)-- (0.85,4.75);
\draw (0.85,4.75)-- (1.38,6.71);
\draw (0.85,4.75)-- (2.08,4.77);
\draw (0.85,4.75)-- (2.08,3.4);
\draw (2.08,3.4)-- (2.08,4.77);
\draw (2.08,3.4)-- (3.14,2.46);
\draw (3.14,2.46)-- (2.08,4.77);
\draw (3.14,2.46)-- (5.41,0.42);
\draw (5.41,0.42)-- (2.08,4.77);
\draw (5.41,0.42)-- (1.38,6.71);
\draw (5.41,0.42)-- (-2.6,0.42);
\draw (-2.6,0.42)-- (1.38,6.71);
\draw (1.38,6.71)-- (2.08,4.77);
\begin{scriptsize}
\fill [color=black] (1.44,2.48) circle (1.5pt);
\fill [color=black] (2,1.68) circle (1.5pt);
\fill [color=black] (0.73,1.68) circle (1.5pt);
\fill [color=ffffff] (0.26,2.48) circle (1.5pt);
\draw (0.26,2.48) circle (1.5pt);
\fill [color=black] (0.89,3.4) circle (1.5pt);
\fill [color=ffffff] (2.22,2.48) circle (1.5pt);
\draw (2.22,2.48) circle (1.5pt);
\fill [color=black] (3.37,1.38) circle (1.5pt);
\fill [color=ffffff] (1.42,0.93) circle (1.5pt);
\draw(1.42,0.93) circle (1.5pt);
\fill [color=black] (-0.46,1.4) circle (1.5pt);
\fill [color=black] (-0.58,2.48) circle (1.5pt);
\fill [color=ffffff] (0.85,4.75) circle (1.5pt);
\draw (0.85,4.75) circle (1.5pt);
\fill [color=black] (2.08,3.4) circle (1.5pt);
\fill [color=black] (3.14,2.46) circle (1.5pt);
\fill [color=ffffff] (5.41,0.42) circle (1.5pt);
\draw (5.41,0.42) circle (1.5pt);
\fill [color=black] (-2.6,0.42) circle (1.5pt);
\fill [color=black] (1.38,6.71) circle (1.5pt);
\fill [color=black] (2.08,4.77) circle (1.5pt);
\end{scriptsize}
\end{tikzpicture}}
\hskip 0.2in
\subfigure[]{ \definecolor{ffffff}{rgb}{1,1,1}
\begin{tikzpicture}[line cap=round,line join=round,>=triangle 45,x=0.32cm,y=0.36cm]
\clip(-2.56,-0.24) rectangle (4.8,5.8);
\draw (0.92,2.86)-- (1.82,1.98);
\draw (0.92,2.86)-- (0.86,1.82);
\draw (0.92,2.86)-- (0.32,3.28);
\draw (0.92,2.86)-- (0.96,4.18);
\draw (0.92,2.86)-- (1.42,3.16);
\draw (1.82,1.98)-- (1.42,3.16);
\draw (1.82,1.98)-- (2.64,1.68);
\draw (1.82,1.98)-- (1.94,1.3);
\draw (1.82,1.98)-- (0.86,1.82);
\draw (0.86,1.82)-- (1.94,1.3);
\draw (0.86,1.82)-- (0.88,1);
\draw (0.86,1.82)-- (0.06,1.8);
\draw (0.86,1.82)-- (0.32,3.28);
\draw (0.32,3.28)-- (0.06,1.8);
\draw (0.32,3.28)-- (-0.6,2.08);
\draw (0.32,3.28)-- (1.1,5.62);
\draw (0.32,3.28)-- (0.96,4.18);
\draw (0.96,4.18)-- (1.1,5.62);
\draw (0.96,4.18)-- (1.88,3.68);
\draw (0.96,4.18)-- (1.42,3.16);
\draw (1.42,3.16)-- (1.88,3.68);
\draw (1.42,3.16)-- (2.64,1.68);
\draw (2.64,1.68)-- (1.88,3.68);
\draw (2.64,1.68)-- (4.66,-0.04);
\draw (2.64,1.68)-- (3.26,0.32);
\draw (2.64,1.68)-- (1.94,1.3);
\draw (1.94,1.3)-- (3.26,0.32);
\draw (1.94,1.3)-- (0.88,1);
\draw (0.88,1)-- (3.26,0.32);
\draw (0.88,1)-- (-0.26,0.66);
\draw (0.88,1)-- (0.06,1.8);
\draw (0.06,1.8)-- (-0.26,0.66);
\draw (0.06,1.8)-- (-0.6,2.08);
\draw (-0.6,2.08)-- (-0.26,0.66);
\draw (-0.6,2.08)-- (-2.42,-0.04);
\draw (-0.6,2.08)-- (1.1,5.62);
\draw (1.1,5.62)-- (-2.42,-0.04);
\draw (1.1,5.62)-- (4.66,-0.04);
\draw (1.1,5.62)-- (1.88,3.68);
\draw (1.88,3.68)-- (4.66,-0.04);
\draw (4.66,-0.04)-- (-2.42,-0.04);
\draw (4.66,-0.04)-- (3.26,0.32);
\draw (3.26,0.32)-- (-2.42,-0.04);
\draw (3.26,0.32)-- (-0.26,0.66);
\draw (-0.26,0.66)-- (-2.42,-0.04);
\begin{scriptsize}
\fill [color=black] (0.92,2.86) circle (1.5pt);
\fill [color=black] (1.82,1.98) circle (1.5pt);
\fill [color=ffffff] (0.86,1.82) circle (1.5pt);
\draw (0.86,1.82) circle (1.5pt);
\fill [color=ffffff] (0.32,3.28) circle (1.5pt);
\draw (0.32,3.28) circle (1.5pt);
\fill [color=black] (0.96,4.18) circle (1.5pt);
\fill [color=black] (1.42,3.16) circle (1.5pt);
\fill [color=ffffff] (2.64,1.68) circle (1.5pt);
\draw (2.64,1.68) circle (1.5pt);
\fill [color=black] (1.94,1.3) circle (1.5pt);
\fill [color=black] (0.88,1) circle (1.5pt);
\fill [color=black] (0.06,1.8) circle (1.5pt);
\fill [color=black] (-0.6,2.08) circle (1.5pt);
\fill [color=ffffff] (1.1,5.62) circle (1.5pt);
\draw (1.1,5.62) circle (1.5pt);
\fill [color=black] (1.88,3.68) circle (1.5pt);
\fill [color=black] (4.66,-0.04) circle (1.5pt);
\fill [color=ffffff] (3.26,0.32) circle (1.5pt);
\draw (3.26,0.32) circle (1.5pt);
\fill [color=black] (-0.26,0.66) circle (1.5pt);
\fill [color=black] (-2.42,-0.04) circle (1.5pt);
\end{scriptsize}
\end{tikzpicture}}
\hskip 0.2in
\subfigure[]{\definecolor{ffffff}{rgb}{1,1,1}
\begin{tikzpicture}[line cap=round,line join=round,>=triangle 45,x=0.38cm,y=0.38cm]
\clip(-1.68,0.16) rectangle (4.8,5.62);
\draw (1.14,2.2)-- (1.66,1.66);
\draw (1.14,2.2)-- (0.72,1.22);
\draw (1.14,2.2)-- (0.3,1.76);
\draw (1.14,2.2)-- (0.74,2.8);
\draw (1.14,2.2)-- (1.38,3.24);
\draw (1.14,2.2)-- (1.94,2.58);
\draw (1.66,1.66)-- (1.94,2.58);
\draw (1.66,1.66)-- (2.9,1.32);
\draw (1.66,1.66)-- (1.86,1.06);
\draw (1.66,1.66)-- (0.72,1.22);
\draw (0.72,1.22)-- (1.86,1.06);
\draw (0.72,1.22)-- (-0.24,0.84);
\draw (0.72,1.22)-- (0.3,1.76);
\draw (0.3,1.76)-- (-0.24,0.84);
\draw (0.3,1.76)-- (-1.5,0.36);
\draw (0.3,1.76)-- (0.74,2.8);
\draw (0.74,2.8)-- (-1.5,0.36);
\draw (0.74,2.8)-- (0.88,3.68);
\draw (0.74,2.8)-- (1.38,3.24);
\draw (1.38,3.24)-- (0.88,3.68);
\draw (1.38,3.24)-- (1.86,3.7);
\draw (1.38,3.24)-- (1.94,2.58);
\draw (1.94,2.58)-- (1.86,3.7);
\draw (1.94,2.58)-- (2.78,2.48);
\draw (1.94,2.58)-- (2.9,1.32);
\draw (2.9,1.32)-- (2.78,2.48);
\draw (2.9,1.32)-- (4.6,0.36);
\draw (2.9,1.32)-- (2.5,0.62);
\draw (2.9,1.32)-- (1.86,1.06);
\draw (1.86,1.06)-- (2.5,0.62);
\draw (1.86,1.06)-- (-0.24,0.84);
\draw (-0.24,0.84)-- (2.5,0.62);
\draw (-0.24,0.84)-- (-1.5,0.36);
\draw (-1.5,0.36)-- (2.5,0.62);
\draw (-1.5,0.36)-- (4.6,0.36);
\draw (-1.5,0.36)-- (1.54,5.48);
\draw (0.88,3.68)-- (-1.5,0.36);
\draw (0.88,3.68)-- (1.54,5.48);
\draw (0.88,3.68)-- (1.86,3.7);
\draw (1.54,5.48)-- (1.86,3.7);
\draw (1.86,3.7)-- (2.78,2.48);
\draw (2.78,2.48)-- (1.54,5.48);
\draw (2.78,2.48)-- (4.6,0.36);
\draw (4.6,0.36)-- (1.54,5.48);
\draw (4.6,0.36)-- (2.5,0.62);
\begin{scriptsize}
\fill [color=ffffff] (1.14,2.2) circle (1.5pt);
\draw (1.14,2.2) circle (1.5pt);
\fill [color=black] (1.66,1.66) circle (1.5pt);
\fill [color=black] (0.72,1.22) circle (1.5pt);
\fill [color=black] (0.3,1.76) circle (1.5pt);
\fill [color=black] (0.74,2.8) circle (1.5pt);
\fill [color=black] (1.38,3.24) circle (1.5pt);
\fill [color=ffffff] (1.94,2.58) circle (1.5pt);
\fill [color=ffffff] (2.9,1.32) circle (1.5pt);
\draw (1.94,2.58) circle (1.5pt);
\draw (2.9,1.32) circle (1.5pt);
\fill [color=black] (1.86,1.06) circle (1.5pt);
\fill [color=black] (-0.24,0.84) circle (1.5pt);
\fill [color=ffffff] (-1.5,0.36) circle (1.5pt);
\draw (-1.5,0.36) circle (1.5pt);
\fill [color=black] (0.88,3.68) circle (1.5pt);
\fill [color=black] (1.86,3.7) circle (1.5pt);
\fill [color=black] (2.78,2.48) circle (1.5pt);
\fill [color=black] (4.6,0.36) circle (1.5pt);
\fill [color=black] (2.5,0.62) circle (1.5pt);
\fill [color=black] (1.54,5.48) circle (1.5pt);
\end{scriptsize}
\end{tikzpicture}}
\caption{All non-isomorphic triangulations of 17 vertices with $\delta=5$. \label{fig20111026-1}}
\end{figure}
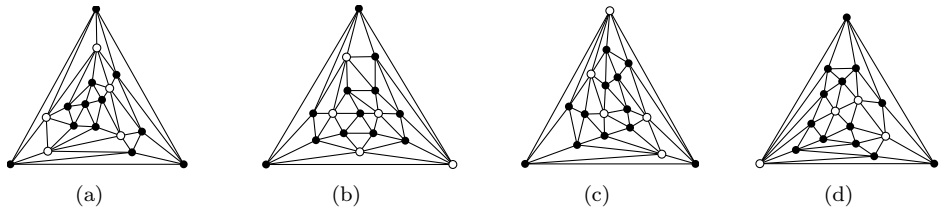

\begin{fact}\label{lem20111026}
There are all together 4 non-isomorphic triangulations of planar graphs on 17 vertices with $\delta=5$ (Figure \ref{fig20111026-1}).
\end{fact}

\begin{fact}\label{lem20111029}
Let $G$ be a triangulation on 18 vertices with minimum degree 5, and let $T$ be the set of vertices whose degree is at least 6 in $G$, then there is no $C_5$ in the subgraph induced by $T$.
\end{fact}

\begin{Lemma}\label{lem20101012}
Let $G$ be a $C_4$-free plane graph. Let $H$ be the subgraph induced by $\Gamma(G)$.

\noindent {\rm (i)} If $G$ is 4-regular, then either $\tau(G)=0$ or $\tau(G)\ge 5$.

\noindent {\rm (ii)} If $G$ is 4-regular and $\tau(G)=5$, then $H$ is a 5-cycle.

\noindent {\rm (iii)} If $\delta(G)\ge 4$, and there is an edge $e=uv \in \Gamma(G)$ such that $d_G(u)=4$, then there is an edge $f=uw$ ($w\ne v$) such that $f\in \Gamma(G)$. 
\end{Lemma}

\begin{figure}[H]
\centering
\begin{tikzpicture}[line cap=round,line join=round,>=triangle 45,x=1.0cm,y=1.0cm]
\clip(-1.38,0.86) rectangle (3.5,4.12);
\draw (-0.94,2.46)-- (-0.22,3.2);
\draw (-0.22,3.2)-- (1.02,3.2);
\draw (1.02,3.2)-- (1.02,1.74);
\draw (1.02,1.74)-- (-0.24,1.74);
\draw (-0.24,1.74)-- (-0.94,2.46);
\draw (1.02,3.2)-- (2.2,3.2);
\draw (2.2,3.2)-- (2.92,2.44);
\draw (2.92,2.44)-- (2.26,1.74);
\draw (2.26,1.74)-- (1.02,1.74);
\draw (1.02,1.74)-- (1.04,0.96);
\draw (1.02,3.2)-- (1.02,4.04);
\draw (0.95,2.32) node[anchor=north west] {$u$};
\draw (0.95,3.2) node[anchor=north west] {$v$};
\draw (-0.38,3.2) node[anchor=north west] {$x_1$};
\draw (1.8,3.2) node[anchor=north west] {$x_2$};
\draw (-0.38,2.32) node[anchor=north west] {$x_3$};
\draw (1.8,2.32) node[anchor=north west] {$x_4$};
\draw (0.28,3.9) node[anchor=north west] {$f_1$};
\draw (1.48,3.9) node[anchor=north west] {$f_2$};
\draw (0.24,1.84) node[anchor=north west] {$f_3$};
\draw (1.44,1.78) node[anchor=north west] {$f_4$};
\draw (-0.22,3.2)-- (-0.34,3.62);
\draw (-0.22,3.2)-- (-0.6,3.36);
\draw (-0.94,2.46)-- (-1.24,2.62);
\draw (-0.94,2.46)-- (-1.24,2.28);
\draw (-0.24,1.74)-- (-0.6,1.54);
\draw (-0.24,1.74)-- (-0.34,1.32);
\draw (2.26,1.74)-- (2.4,1.26);
\draw (2.26,1.74)-- (2.6,1.52);
\draw (2.92,2.44)-- (3.34,2.12);
\draw (2.92,2.44)-- (3.38,2.54);
\draw (2.2,3.2)-- (2.72,3.28);
\draw (2.2,3.2)-- (2.48,3.7);
\begin{scriptsize}
\fill [color=black] (-0.94,2.46) circle (1.5pt);
\fill [color=black] (-0.22,3.2) circle (1.5pt);
\fill [color=black] (1.02,3.2) circle (1.5pt);
\fill [color=black] (1.02,1.74) circle (1.5pt);
\fill [color=black] (-0.24,1.74) circle (1.5pt);
\fill [color=black] (2.2,3.2) circle (1.5pt);
\fill [color=black] (2.92,2.44) circle (1.5pt);
\fill [color=black] (2.26,1.74) circle (1.5pt);
\end{scriptsize}
\end{tikzpicture}
\caption{$\tau(G)\ge 5$}\label{fig20101012}
\end{figure}
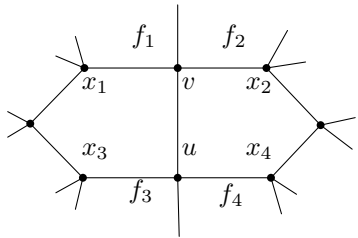

 {\bf Proof.}
Suppose that $G$ is 4-regular and $\tau(G)\ne 0$, let $uv$ be an arbitrary edge in $\Gamma(G)$, then there exists two faces $g_1,g_2\in F(G)-F_3(G)$ such that $g_1, g_2$ have an edge $uv$ in common (see Figure \ref{fig20101012}).  Since $G$ is 4-regular, $f_1\ne f_2$ and $f_3\ne f_4$. Furthermore, since $G$ is $C_4$-free, at least one of $f_1$ and $f_2$ are non-triangles, so at least one of $vx_1$ and $vx_2$ are not covered by triangles. For the same reason as above,  we see that at least one of $ux_3$ and $ux_4$ are not covered by triangles. This implies $H$ has at least 4 vertices and that both $d_H(u)\ge 2$ and $d_H(v)\ge 2$, and hence $\delta(H)\ge 2$. If $\tau(G)= 4$, then $H$ will contain a 4-cycle, a contradiction. So we have $\tau(G)\ge 5$. This complete the proof of (i).

By the above arguments, we see that $\delta(H)\ge 2$ and it is obvious that (ii) and (iii) holds.

\begin{Lemma}\label{lem20111029-1}
Let $G$ be a 4-regular $C_4$-free plane graph with $\tau(G)=5, f_6\le 1, f_7=f_8=\cdots=0$, then $\Gamma(G)$ does not induce a separating 5-cycle in $G$.
\end{Lemma}

{\bf Proof.} Let $H$ be the subgraph induced by $\Gamma(G)$. By the proof of Lemma \ref{lem20101012},  we know that $\delta(H)\ge 2$. Since $G$ is $C_4$-free, $H$ must be a 5-cycle $C$ in $G$. Suppose on the contrary that $C$ is a separating 5-cycle.

{\bf Claim.} For each vertex $v$ in $V(C)$, the two edges which are adjacent to $v$ and which are not on $C$ must be  either outside or inside $C$, but not both.

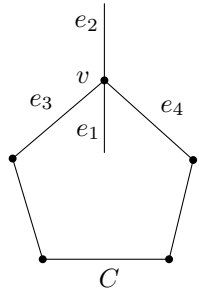
\begin{figure}[H]
\centering
\begin{tikzpicture}[line cap=round,line join=round,>=triangle 45,x=1.0cm,y=1.0cm]
\clip(-0.6,0.62) rectangle (2.2,4.74);
\draw (0.82,3.62)-- (-0.4,2.58);
\draw (0.82,3.62)-- (2,2.56);
\draw (-0.4,2.58)-- (0,1.24);
\draw (0,1.24)-- (1.68,1.24);
\draw (2,2.56)-- (1.68,1.24);
\draw (0.82,3.62)-- (0.82,4.64);
\draw (0.82,3.62)-- (0.82,2.66);
\draw (0.32,3.86) node[anchor=north west] {$v$};
\draw (0.32,3.16) node[anchor=north west] {$e_1$};
\draw (0.3,4.6) node[anchor=north west] {$e_2$};
\draw (-0.3,3.56) node[anchor=north west] {$e_3$};
\draw (1.44,3.48) node[anchor=north west] {$e_4$};
\draw (0.62,1.24) node[anchor=north west] {$C$};
\begin{scriptsize}
\fill [color=black] (0.82,3.62) circle (1.5pt);
\fill [color=black] (-0.4,2.58) circle (1.5pt);
\fill [color=black] (2,2.56) circle (1.5pt);
\fill [color=black] (0,1.24) circle (1.5pt);
\fill [color=black] (1.68,1.24) circle (1.5pt);
\end{scriptsize}
\end{tikzpicture}
\caption{A forbidden structure\label{fig20111027}}
\end{figure}

{\bf Proof of Claim.} Let $e_1, e_2$ be the  two edges which are incident with $v$ and which are not on $C$. Suppose, without loss of generality, that $e_1$ is  inside $C$,  and  $e_2$ is  outside $C$ (Figure \ref{fig20111027}). Since $e_1$ is covered by a triangle in $G$, and $e_3,e_4$ are not covered by any triangle, this is impossible since $G$ is 4-regular and $C_4$-free. $\qed$

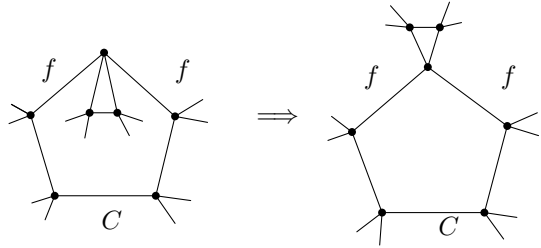
\begin{figure}[H]
\centering
\begin{tikzpicture}[line cap=round,line join=round,>=triangle 45,x=0.8cm,y=0.8cm]
\clip(-0.86,0.3) rectangle (8.08,4.56);
\draw (0.82,3.62)-- (-0.4,2.58);
\draw (0.82,3.62)-- (2,2.56);
\draw (-0.4,2.58)-- (0,1.24);
\draw (0,1.24)-- (1.68,1.24);
\draw (2,2.56)-- (1.68,1.24);
\draw (0.82,3.62)-- (1.04,2.62);
\draw (0.82,3.62)-- (0.58,2.62);
\draw (0.62,1.14) node[anchor=north west] {$C$};
\draw (0.58,2.62)-- (1.04,2.62);
\draw (0.58,2.62)-- (0.18,2.48);
\draw (0.58,2.62)-- (0.5,2.2);
\draw (1.04,2.62)-- (1.24,2.24);
\draw (1.04,2.62)-- (1.46,2.48);
\draw (-0.4,2.58)-- (-0.72,2.76);
\draw (-0.4,2.58)-- (-0.72,2.76);
\draw (-0.4,2.58)-- (-0.76,2.44);
\draw (0,1.24)-- (-0.38,1.12);
\draw (0,1.24)-- (-0.16,0.84);
\draw (1.68,1.24)-- (2,0.78);
\draw (1.68,1.24)-- (2.26,1.16);
\draw (2,2.56)-- (2.46,2.84);
\draw (2,2.56)-- (2.54,2.46);
\draw (1.84,3.66) node[anchor=north west] {$ f $};
\draw (3.2,2.76) node[anchor=north west] {$ \Longrightarrow $};
\draw (6.2,3.38)-- (4.94,2.3);
\draw (6.2,3.38)-- (7.52,2.4);
\draw (4.94,2.3)-- (5.44,0.96);
\draw (7.52,2.4)-- (7.14,0.96);
\draw (5.44,0.96)-- (7.14,0.96);
\draw (6.2,3.38)-- (5.9,4.04);
\draw (6.2,3.38)-- (6.4,4.04);
\draw (5.9,4.04)-- (6.4,4.04);
\draw (6.4,4.04)-- (6.56,4.46);
\draw (6.4,4.04)-- (6.76,4.12);
\draw (5.9,4.04)-- (5.58,4.4);
\draw (5.9,4.04)-- (5.5,4.08);
\draw (4.94,2.3)-- (4.6,2.56);
\draw (4.94,2.3)-- (4.54,2.16);
\draw (7.52,2.4)-- (8.02,2.62);
\draw (7.52,2.4)-- (8.04,2.24);
\draw (7.14,0.96)-- (7.46,0.52);
\draw (7.14,0.96)-- (7.68,0.88);
\draw (5.44,0.96)-- (5.02,0.64);
\draw (5.44,0.96)-- (5.4,0.42);
\draw (7.26,3.56) node[anchor=north west] {$f$};
\draw (-0.38,3.68) node[anchor=north west] {$f$};
\draw (4.98,3.56) node[anchor=north west] {$f$};
\draw (6.22,1.04) node[anchor=north west] {$C$};
\begin{scriptsize}
\fill [color=black] (0.82,3.62) circle (1.5pt);
\fill [color=black] (-0.4,2.58) circle (1.5pt);
\fill [color=black] (2,2.56) circle (1.5pt);
\fill [color=black] (0,1.24) circle (1.5pt);
\fill [color=black] (1.68,1.24) circle (1.5pt);
\fill [color=black] (1.04,2.62) circle (1.5pt);
\fill [color=black] (0.58,2.62) circle (1.5pt);
\fill [color=black] (6.2,3.38) circle (1.5pt);
\fill [color=black] (4.94,2.3) circle (1.5pt);
\fill [color=black] (7.52,2.4) circle (1.5pt);
\fill [color=black] (5.44,0.96) circle (1.5pt);
\fill [color=black] (7.14,0.96) circle (1.5pt);
\fill [color=black] (5.9,4.04) circle (1.5pt);
\fill [color=black] (6.4,4.04) circle (1.5pt);
\end{scriptsize}
\end{tikzpicture}
\caption{Re-embedding  of $G$ }\label{fig20111027-1}
\end{figure}

Let $t_1$ and $t_2$ be the number of edges that are incident with $V(C)$ and belong to  the inside of $C$ and the outside of $C$ respectively. Since $G$ is 4-regular, by Claim we see that $t_1+t_2=10$, and both $t_1$ and $t_2$ are even integers.
Without loss of generality, we assume that $t_1\le 4$.  Suppose first that $t_1=2$ (Figure \ref{fig20111027-1}). In this case we can embed all the vertices of inside $C$ to outside $C$, this contradicts that $C$ is a separating cycle.

\begin{figure}[H]
\centering
\begin{tikzpicture}[line cap=round,line join=round,>=triangle 45,x=1.0cm,y=1.0cm]
\clip(-0.86,0.3) rectangle (8.16,4.56);
\draw (0.9,3.6)-- (-0.4,2.58);
\draw (0.9,3.6)-- (2.26,2.56);
\draw (-0.4,2.58)-- (-0.04,0.96);
\draw (-0.04,0.96)-- (1.76,0.94);
\draw (2.26,2.56)-- (1.76,0.94);
\draw (0.9,3.6)-- (1.34,2.96);
\draw (0.9,3.6)-- (1.04,2.8);
\draw (0.62,0.94) node[anchor=north west] {$C$};
\draw (1.04,2.8)-- (1.34,2.96);
\draw (1.04,2.8)-- (0.96,2.44);
\draw (1.04,2.8)-- (1.14,2.46);
\draw (1.34,2.96)-- (1.42,2.5);
\draw (1.34,2.96)-- (1.62,2.62);
\draw (-0.4,2.58)-- (0.52,2.06);
\draw (-0.4,2.58)-- (0.52,2.06);
\draw (-0.4,2.58)-- (0.28,1.72);
\draw (-0.04,0.96)-- (-0.38,1.12);
\draw (-0.04,0.96)-- (-0.28,0.7);
\draw (1.76,0.94)-- (2.06,0.74);
\draw (1.76,0.94)-- (2.12,1.08);
\draw (2.26,2.56)-- (2.68,2.76);
\draw (2.26,2.56)-- (2.68,2.42);
\draw (1.92,3.82) node[anchor=north west] {$ f $};
\draw (3.26,1.98) node[anchor=north west] {$ \Longrightarrow $};
\draw (6.2,3.38)-- (5.02,2.38);
\draw (6.2,3.38)-- (7.52,2.4);
\draw (5.02,2.38)-- (5.44,0.96);
\draw (7.52,2.4)-- (7.14,0.96);
\draw (5.44,0.96)-- (7.14,0.96);
\draw (6.2,3.38)-- (5.9,4.04);
\draw (6.2,3.38)-- (6.4,4.04);
\draw (5.9,4.04)-- (6.4,4.04);
\draw (6.4,4.04)-- (6.56,4.46);
\draw (6.4,4.04)-- (6.76,4.12);
\draw (5.9,4.04)-- (5.58,4.4);
\draw (5.9,4.04)-- (5.5,4.08);
\draw (5.02,2.38)-- (4.44,3.1);
\draw (5.02,2.38)-- (4.18,2.74);
\draw (7.52,2.4)-- (8.02,2.62);
\draw (7.52,2.4)-- (8.04,2.24);
\draw (7.14,0.96)-- (7.46,0.52);
\draw (7.14,0.96)-- (7.68,0.88);
\draw (5.44,0.96)-- (5.02,0.64);
\draw (5.44,0.96)-- (5.4,0.42);
\draw (7.26,3.56) node[anchor=north west] {$f$};
\draw (-0.38,3.68) node[anchor=north west] {$f$};
\draw (4.98,3.56) node[anchor=north west] {$f$};
\draw (6,0.94) node[anchor=north west] {$C$};
\draw (0.28,1.72)-- (0.52,2.06);
\draw (0.52,2.06)-- (0.84,1.84);
\draw (0.52,2.06)-- (0.9,2.12);
\draw (0.28,1.72)-- (0.46,1.38);
\draw (0.28,1.72)-- (0.64,1.58);
\draw (4.22,2.72)-- (4.44,3.1);
\draw (4.44,3.1)-- (4.4,3.48);
\draw (4.44,3.1)-- (4.7,3.38);
\draw (4.22,2.72)-- (3.88,2.9);
\draw (4.22,2.72)-- (3.9,2.56);
\draw (-0.92,2.22) node[anchor=north west] {$f$};
\draw (4.72,2.08) node[anchor=north west] {$f$};
\draw (5.22,2.52) node[anchor=north west] {$u$};
\draw (6,3.16) node[anchor=north west] {$v$};
\draw (-0.94,2.85) node[anchor=north west] {$u$};
\draw (0.8,4.06) node[anchor=north west] {$v$};
\begin{scriptsize}
\fill [color=black] (0.9,3.6) circle (1.5pt);
\fill [color=black] (-0.4,2.58) circle (1.5pt);
\fill [color=black] (2.26,2.56) circle (1.5pt);
\fill [color=black] (-0.04,0.96) circle (1.5pt);
\fill [color=black] (1.76,0.94) circle (1.5pt);
\fill [color=black] (1.34,2.96) circle (1.5pt);
\fill [color=black] (1.04,2.8) circle (1.5pt);
\fill [color=black] (0.28,1.72) circle (1.5pt);
\fill [color=black] (6.2,3.38) circle (1.5pt);
\fill [color=black] (5.02,2.38) circle (1.5pt);
\fill [color=black] (7.52,2.4) circle (1.5pt);
\fill [color=black] (5.44,0.96) circle (1.5pt);
\fill [color=black] (7.14,0.96) circle (1.5pt);
\fill [color=black] (5.9,4.04) circle (1.5pt);
\fill [color=black] (6.4,4.04) circle (1.5pt);
\fill [color=black] (4.44,3.1) circle (1.5pt);
\fill [color=black] (4.18,2.74) circle (1.5pt);
\fill [color=black] (0.52,2.06) circle (1.5pt);
\fill [color=black] (0.28,1.72) circle (1.5pt);
\fill [color=black] (4.44,3.1) circle (1.5pt);
\end{scriptsize}
\end{tikzpicture}
\caption{Another re-embedding of $G$\label{fig20111027-2}}
\end{figure}
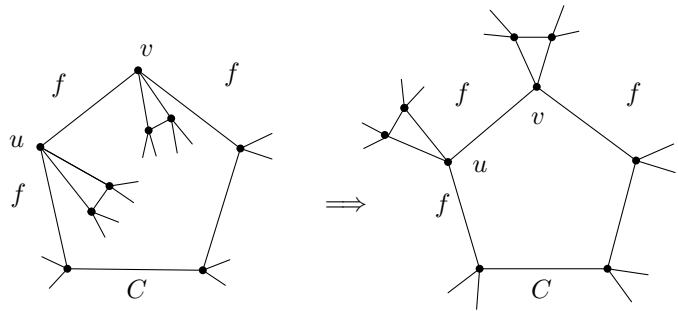

Suppose next that $t_1=4$, let  $u,v$ be the two vertices on $C$ such that the edges which are not on $C$ and are incident with $u,v$  are inside $C$. If $u, v$ are adjacent on $C$, then we can also re-embed all the vertices inside $C$ to outside $C$ (Figure \ref{fig20111027-2}), this contradicts again that $C$ is a separating cycle. 

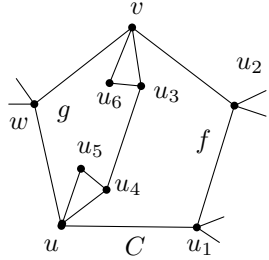
\begin{figure}[H]
\centering
\begin{tikzpicture}[line cap=round,line join=round,>=triangle 45,x=1.0cm,y=1.0cm]
\clip(-1.06,0.18) rectangle (2.88,4.14);
\draw (0.9,3.6)-- (-0.4,2.58);
\draw (0.9,3.6)-- (2.26,2.56);
\draw (-0.4,2.58)-- (-0.04,0.96);
\draw (-0.04,0.96)-- (1.76,0.94);
\draw (2.26,2.56)-- (1.76,0.94);
\draw (0.9,3.6)-- (1.02,2.82);
\draw (0.9,3.6)-- (0.6,2.86);
\draw (0.68,0.94) node[anchor=north west] {$C$};
\draw (0.6,2.86)-- (1.02,2.82);
\draw (1.76,0.94)-- (2.06,0.74);
\draw (1.76,0.94)-- (2.12,1.08);
\draw (2.26,2.56)-- (2.68,2.76);
\draw (2.26,2.56)-- (2.68,2.42);
\draw (1.62,2.36) node[anchor=north west] {$ f $};
\draw (-0.4,0.9) node[anchor=north west] {$u$};
\draw (0.75,4.06) node[anchor=north west] {$v$};
\draw (-0.04,0.96)-- (0.22,1.72);
\draw (-0.03,0.99)-- (0.56,1.44);
\draw (0.22,1.72)-- (0.56,1.44);
\draw (-0.4,2.58)-- (-0.64,2.92);
\draw (-0.4,2.58)-- (-0.74,2.58);
\draw (1.02,2.82)--(0.56,1.44);
\draw (1.5,0.9) node[anchor=north west] {$u_1$};
\draw (2.16,3.32) node[anchor=north west] {$u_2$};
\draw (1.06,2.94) node[anchor=north west] {$u_3$};
\draw (0.54,1.7) node[anchor=north west] {$u_4$};
\draw (-0.22,2.68) node[anchor=north west] {$g$};
\draw (0.04,2.2) node[anchor=north west] {$u_5$};
\draw (0.3,2.86) node[anchor=north west] {$u_6$};
\draw (-0.86,2.52) node[anchor=north west] {$w$};
\begin{scriptsize}
\fill [color=black] (0.9,3.6) circle (1.5pt);
\fill [color=black] (-0.4,2.58) circle (1.5pt);
\fill [color=black] (2.26,2.56) circle (1.5pt);
\fill [color=black] (-0.04,0.96) circle (1.5pt);
\fill [color=black] (1.76,0.94) circle (1.5pt);
\fill [color=black] (1.02,2.82) circle (1.5pt);
\fill [color=black] (0.6,2.86) circle (1.5pt);
\fill [color=black] (0.22,1.72) circle (1.5pt);
\fill [color=black] (-0.03,0.99) circle (1.5pt);
\fill [color=black] (0.56,1.44) circle (1.5pt);
\end{scriptsize}
\end{tikzpicture}
\caption{An forbidden structure of $G$.\label{fig20111027-3}}
\end{figure}

So we assume that $u,v$ are not adjacent on $C$ (Figure \ref{fig20111027-3}). In this case, consider the face $f$ inside $C$ which  are incident with $u_1u_2$. Since $G$ is $C_4$-free, $u_3\ne u_4$. Hence $u,u_1,u_2,v,u_3,u_4$ are all on the boundary of $f$, which implies that the length of $f$ is at least $6$, so 
$f$ must be a 6-face because of the initial hypothesis that for each $k\ge 7$, $f_k(G)=0$. This implies that $u_3$ and $u_4$ are adjacent. Consider the face $g$ which is incident with $w$ and inside $C$, since $f_6\le 1$ and $\{w,u,u_5,u_6,v\}$ is on the boundary of $g$, this implies that $u_5$ and $u_6$ are adjacent, but now $u_3u_4u_5u_6u_3$ is a $C_4$ in $G$, which contradicts the initial hypothesis that $G$ is $C_4$-free. \qed

\subsection{Proof of Theorem \ref{thm20101011}}
{\bf Proof.} (i) By Corollary \ref{coro20101010}, it suffices to show that for each $n\in A$, there is a $C_4$-free planar graph $G$ of order $n$ which has minimum degree 4. 

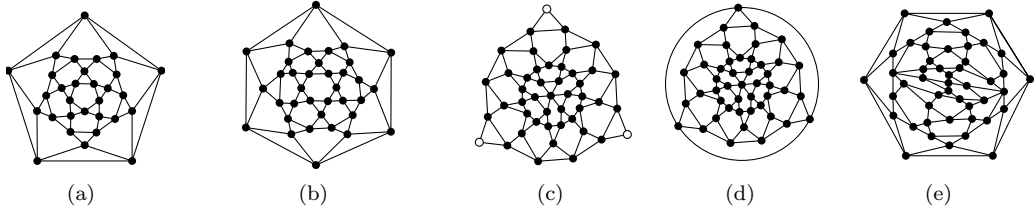
\begin{figure}[H]
\centering 
\subfigure[]{ 
\begin{tikzpicture}[line cap=round,line join=round,>=triangle 45,x=0.21cm,y=0.21cm]
\clip(-2.98,-3.5) rectangle (7.08,6.44);
\draw (1.42,1.9)-- (2.58,1.88);
\draw (1.42,1.9)-- (1.06,0.84);
\draw (1.06,0.84)-- (2,0.14);
\draw (2,0.14)-- (2.92,0.82);
\draw (2.58,1.88)-- (2.92,0.82);
\draw (1.42,1.9)-- (1.98,2.68);
\draw (1.98,2.68)-- (2.58,1.88);
\draw (2.58,1.88)-- (3.46,1.6);
\draw (3.46,1.6)-- (2.92,0.82);
\draw (2.92,0.82)-- (2.92,-0.18);
\draw (2.92,-0.18)-- (2,0.14);
\draw (1.06,0.84)-- (1.12,-0.14);
\draw (1.12,-0.14)-- (2,0.14);
\draw (1.06,0.84)-- (0.54,1.62);
\draw (0.54,1.62)-- (1.42,1.9);
\draw (1.98,2.68)-- (2.76,3.46);
\draw (2.76,3.46)-- (4.02,2.48);
\draw (3.46,1.6)-- (4.02,2.48);
\draw (3.46,1.6)-- (4.44,1.16);
\draw (2.92,-0.18)-- (3.92,-0.34);
\draw (4.44,1.16)-- (3.92,-0.34);
\draw (2.92,-0.18)-- (2.74,-1.2);
\draw (1.12,-0.14)-- (1.1,-1.16);
\draw (1.1,-1.16)-- (2.74,-1.2);
\draw (1.12,-0.14)-- (-0.04,-0.24);
\draw (0.54,1.62)-- (-0.44,1.14);
\draw (-0.44,1.14)-- (-0.04,-0.24);
\draw (0.54,1.62)-- (-0.06,2.44);
\draw (-0.06,2.44)-- (1.36,3.46);
\draw (1.98,2.68)-- (1.36,3.46);
\draw (1.36,3.46)-- (2.76,3.46);
\draw (4.02,2.48)-- (4.44,1.16);
\draw (3.92,-0.34)-- (2.74,-1.2);
\draw (1.1,-1.16)-- (-0.04,-0.24);
\draw (-0.44,1.14)-- (-0.06,2.44);
\draw (-0.06,2.44)-- (0.18,3.68);
\draw (0.18,3.68)-- (1.36,3.46);
\draw (2.76,3.46)-- (3.84,3.66);
\draw (3.84,3.66)-- (4.02,2.48);
\draw (4.44,1.16)-- (5,0.16);
\draw (3.92,-0.34)-- (5,0.16);
\draw (1.1,-1.16)-- (2,-2);
\draw (2,-2)-- (2.74,-1.2);
\draw (-0.44,1.14)-- (-1.02,0.16);
\draw (-1.02,0.16)-- (-0.04,-0.24);
\draw (2,6.22)-- (0.18,3.68);
\draw (2,6.22)-- (3.84,3.66);
\draw (6.86,2.72)-- (3.84,3.66);
\draw (6.86,2.72)-- (5,0.16);
\draw (5,-3)-- (2,-2);
\draw (5,-3)-- (5,0.16);
\draw (-1,-3)-- (-1.02,0.16);
\draw (-1,-3)-- (2,-2);
\draw (-1.02,0.16)-- (-2.84,2.72);
\draw (-2.84,2.72)-- (0.18,3.68);
\draw (2,6.22)-- (6.86,2.72);
\draw (6.86,2.72)-- (5,-3);
\draw (5,-3)-- (-1,-3);
\draw (-2.84,2.72)-- (-1,-3);
\draw (-2.84,2.72)-- (2,6.22);
\begin{scriptsize}
\fill [color=black] (2.58,1.88) circle (1.5pt);
\fill [color=black] (1.42,1.9) circle (1.5pt);
\fill [color=black] (1.06,0.84) circle (1.5pt);
\fill [color=black] (2,0.14) circle (1.5pt);
\fill [color=black] (2.92,0.82) circle (1.5pt);
\fill [color=black] (1.98,2.68) circle (1.5pt);
\fill [color=black] (3.46,1.6) circle (1.5pt);
\fill [color=black] (2.92,-0.18) circle (1.5pt);
\fill [color=black] (1.12,-0.14) circle (1.5pt);
\fill [color=black] (0.54,1.62) circle (1.5pt);
\fill [color=black] (2.76,3.46) circle (1.5pt);
\fill [color=black] (4.02,2.48) circle (1.5pt);
\fill [color=black] (4.44,1.16) circle (1.5pt);
\fill [color=black] (3.92,-0.34) circle (1.5pt);
\fill [color=black] (2.74,-1.2) circle (1.5pt);
\fill [color=black] (1.1,-1.16) circle (1.5pt);
\fill [color=black] (-0.04,-0.24) circle (1.5pt);
\fill [color=black] (-0.44,1.14) circle (1.5pt);
\fill [color=black] (-0.06,2.44) circle (1.5pt);
\fill [color=black] (1.36,3.46) circle (1.5pt);
\fill [color=black] (0.18,3.68) circle (1.5pt);
\fill [color=black] (3.84,3.66) circle (1.5pt);
\fill [color=black] (5,0.16) circle (1.5pt);
\fill [color=black] (2,-2) circle (1.5pt);
\fill [color=black] (-1.02,0.16) circle (1.5pt);
\fill [color=black] (2,6.22) circle (1.5pt);
\fill [color=black] (6.86,2.72) circle (1.5pt);
\fill [color=black] (5,-3) circle (1.5pt);
\fill [color=black] (-1,-3) circle (1.5pt);
\fill [color=black] (-2.84,2.72) circle (1.5pt);
\end{scriptsize}
\end{tikzpicture}
} 
\hspace{0.2in} 
\subfigure[]{ 
\begin{tikzpicture}[line cap=round,line join=round,>=triangle 45,x=0.2cm,y=0.2cm]
\clip(-1.33,-5.0) rectangle (9.25,6.3);
\draw (3.36,1.3)-- (4.66,1.32);
\draw (4.66,1.32)-- (5.18,0.4);
\draw (5.18,0.4)-- (4.61,-0.51);
\draw (4.61,-0.51)-- (3.44,-0.47);
\draw (3.44,-0.47)-- (2.84,0.36);
\draw (2.84,0.36)-- (3.36,1.3);
\draw (4.66,1.32)-- (5.64,1.38);
\draw (5.64,1.38)-- (5.18,0.4);
\draw (4.61,-0.51)-- (4.18,-1.4);
\draw (4.18,-1.4)-- (3.44,-0.47);
\draw (2.84,0.36)-- (2.28,1.38);
\draw (2.28,1.38)-- (3.36,1.3);
\draw (5.64,1.38)-- (6.06,2.08);
\draw (6.06,2.08)-- (4.66,2.84);
\draw (2.28,1.38)-- (2,2.12);
\draw (2,2.12)-- (3.3,2.8);
\draw (2.28,1.38)-- (1.48,1.3);
\draw (1.48,1.3)-- (1.4,-0.62);
\draw (1.96,-1.68)-- (3.4,-2.36);
\draw (3.4,-2.36)-- (4.18,-1.4);
\draw (4.18,-1.4)-- (4.86,-2.28);
\draw (4.86,-2.28)-- (6.26,-1.44);
\draw (6.82,-0.44)-- (6.68,0.9);
\draw (6.68,0.9)-- (5.64,1.38);
\draw (4.66,2.84)-- (3.3,2.8);
\draw (6.06,2.08)-- (6.68,0.9);
\draw (6.82,-0.44)-- (6.26,-1.44);
\draw (4.86,-2.28)-- (3.4,-2.36);
\draw (1.96,-1.68)-- (1.4,-0.62);
\draw (1.48,1.3)-- (2,2.12);
\draw (2,2.12)-- (2.06,3.34);
\draw (2.06,3.34)-- (3.3,2.8);
\draw (6.68,0.9)-- (7.68,0.3);
\draw (7.68,0.3)-- (6.82,-0.44);
\draw (3.4,-2.36)-- (2.18,-2.7);
\draw (2.18,-2.7)-- (1.96,-1.68);
\draw (-0.8,2.8)-- (2.06,3.34);
\draw (-0.84,-2.26)-- (2.18,-2.7);
\draw (3.82,-4.76)-- (2.18,-2.7);
\draw (8.8,-2.54)-- (7.68,0.3);
\draw (8.82,2.72)-- (7.68,0.3);
\draw (3.82,5.9)-- (2.06,3.34);
\draw (3.82,5.9)-- (8.82,2.72);
\draw (8.82,2.72)-- (8.8,-2.54);
\draw (8.8,-2.54)-- (3.82,-4.76);
\draw (3.82,-4.76)-- (-0.84,-2.26);
\draw (-0.84,-2.26)-- (-0.8,2.8);
\draw (-0.8,2.8)-- (3.82,5.9);
\draw (5.71,3.14)-- (4.66,2.84);
\draw (5.71,3.14)-- (6.06,2.08);
\draw (0.3,0.33)-- (1.48,1.3);
\draw (0.3,0.33)-- (1.4,-0.62);
\draw (5.91,-2.71)-- (4.86,-2.28);
\draw (5.91,-2.71)-- (6.26,-1.44);
\draw (4.01,2.03)-- (3.36,1.3);
\draw (4.01,2.03)-- (4.66,1.32);
\draw (5.61,-0.47)-- (5.18,0.4);
\draw (5.61,-0.47)-- (4.61,-0.51);
\draw (2.37,-0.67)-- (2.84,0.36);
\draw (2.37,-0.67)-- (3.44,-0.47);
\draw (5.71,3.14)-- (3.82,5.9);
\draw (5.71,3.14)-- (8.82,2.72);
\draw (4.01,2.03)-- (3.3,2.8);
\draw (4.01,2.03)-- (4.66,2.84);
\draw (0.3,0.33)-- (-0.8,2.8);
\draw (0.3,0.33)-- (-0.84,-2.26);
\draw (2.37,-0.67)-- (1.4,-0.62);
\draw (2.37,-0.67)-- (1.96,-1.68);
\draw (5.61,-0.47)-- (6.26,-1.44);
\draw (5.61,-0.47)-- (6.82,-0.44);
\draw (5.91,-2.71)-- (3.82,-4.76);
\draw (5.91,-2.71)-- (8.8,-2.54);
\begin{scriptsize}
\fill [color=black] (3.36,1.3) circle (1.5pt);
\fill [color=black] (4.66,1.32) circle (1.5pt);
\fill [color=black] (5.18,0.4) circle (1.5pt);
\fill [color=black] (4.61,-0.51) circle (1.5pt);
\fill [color=black] (3.44,-0.47) circle (1.5pt);
\fill [color=black] (2.84,0.36) circle (1.5pt);
\fill [color=black] (5.64,1.38) circle (1.5pt);
\fill [color=black] (4.18,-1.4) circle (1.5pt);
\fill [color=black] (2.28,1.38) circle (1.5pt);
\fill [color=black] (6.06,2.08) circle (1.5pt);
\fill [color=black] (4.66,2.84) circle (1.5pt);
\fill [color=black] (3.3,2.8) circle (1.5pt);
\fill [color=black] (2,2.12) circle (1.5pt);
\fill [color=black] (1.48,1.3) circle (1.5pt);
\fill [color=black] (1.4,-0.62) circle (1.5pt);
\fill [color=black] (1.96,-1.68) circle (1.5pt);
\fill [color=black] (3.4,-2.36) circle (1.5pt);
\fill [color=black] (4.86,-2.28) circle (1.5pt);
\fill [color=black] (6.26,-1.44) circle (1.5pt);
\fill [color=black] (6.82,-0.44) circle (1.5pt);
\fill [color=black] (6.68,0.9) circle (1.5pt);
\fill [color=black] (2.06,3.34) circle (1.5pt);
\fill [color=black] (7.68,0.3) circle (1.5pt);
\fill [color=black] (2.18,-2.7) circle (1.5pt);
\fill [color=black] (-0.8,2.8) circle (1.5pt);
\fill [color=black] (-0.84,-2.26) circle (1.5pt);
\fill [color=black] (3.82,-4.76) circle (1.5pt);
\fill [color=black] (8.8,-2.54) circle (1.5pt);
\fill [color=black] (8.82,2.72) circle (1.5pt);
\fill [color=black] (3.82,5.9) circle (1.5pt);
\fill [color=black] (5.71,3.14) circle (1.5pt);
\fill [color=black] (0.3,0.33) circle (1.5pt);
\fill [color=black] (5.91,-2.71) circle (1.5pt);
\fill [color=black] (4.01,2.03) circle (1.5pt);
\fill [color=black] (5.61,-0.47) circle (1.5pt);
\fill [color=black] (2.37,-0.67) circle (1.5pt);
\end{scriptsize}
\end{tikzpicture}
} 
\hspace{0.2in} 
\subfigure[]{\definecolor{ffffff}{rgb}{1,1,1}
\begin{tikzpicture}[line cap=round,line join=round,>=triangle 45,x=0.09cm,y=0.09cm]
\clip(-7.91,-10.96) rectangle (16.87,14);
\draw (0.89,9.3)-- (6.01,9.61);
\draw (0.89,9.3)-- (1.33,6.41);
\draw (1.33,6.41)-- (3.61,4.31);
\draw (3.61,4.31)-- (6.19,6.59);
\draw (6.01,9.61)-- (6.19,6.59);
\draw (6.19,6.59)-- (10.56,7.38);
\draw (10.56,7.38)-- (6.01,9.61);
\draw (0.89,9.3)-- (-3.4,6.5);
\draw (-3.4,6.5)-- (1.33,6.41);
\draw (1.33,6.41)-- (1.86,3.7);
\draw (1.86,3.7)-- (3.61,4.31);
\draw (3.61,4.31)-- (5.49,4.09);
\draw (6.19,6.59)-- (5.49,4.09);
\draw (10.56,7.38)-- (10.04,3.57);
\draw (10.56,7.38)-- (13.45,3);
\draw (10.04,3.57)-- (13.45,3);
\draw (10.04,3.57)-- (7.11,3.13);
\draw (5.49,4.09)-- (7.11,3.13);
\draw (5.49,4.09)-- (5.09,1.99);
\draw (5.09,1.99)-- (7.11,3.13);
\draw (2.42,1.77)-- (1.86,3.7);
\draw (2.42,1.77)-- (5.09,1.99);
\draw (13.45,3)-- (11,-1.11);
\draw (13.45,3)-- (13.71,-2.08);
\draw (11,-1.11)-- (13.71,-2.08);
\draw (13.71,-2.08)-- (11.4,-6.67);
\draw (7.11,3.13)-- (8.2,1.6);
\draw (10.04,3.57)-- (8.2,1.6);
\draw (8.2,1.6)-- (8.55,-0.28);
\draw (8.2,1.6)-- (6.41,0.07);
\draw (6.41,0.07)-- (8.55,-0.28);
\draw (8.55,-0.28)-- (11,-1.11);
\draw (8.55,-0.28)-- (8.16,-2.03);
\draw (11,-1.11)-- (8.16,-2.03);
\draw (8.16,-2.03)-- (7.06,-3.61);
\draw (8.16,-2.03)-- (8.73,-5.45);
\draw (7.06,-3.61)-- (8.73,-5.45);
\draw (11.4,-6.67)-- (8.73,-5.45);
\draw (8.73,-5.45)-- (7.19,-9.47);
\draw (11.4,-6.67)-- (7.19,-9.47);
\draw (1.86,3.7)-- (0.41,2.47);
\draw (0.41,2.47)-- (2.42,1.77);
\draw (-3.4,6.5)-- (-2.35,3.17);
\draw (-2.35,3.17)-- (0.41,2.47);
\draw (6.41,0.07)-- (5.44,-2.16);
\draw (5.44,-2.16)-- (7.06,-3.61);
\draw (-2.35,3.17)-- (-0.29,0.77);
\draw (0.41,2.47)-- (-0.29,0.77);
\draw (-0.29,0.77)-- (1.55,-0.37);
\draw (1.55,-0.37)-- (3.34,-2.56);
\draw (-3.4,6.5)-- (-5.71,1.86);
\draw (-2.35,3.17)-- (-5.71,1.86);
\draw (-0.29,0.77)-- (-0.33,-1.11);
\draw (-0.33,-1.11)-- (1.55,-0.37);
\draw (-5.71,1.86)-- (-3,-1.64);
\draw (-3,-1.64)-- (-0.33,-1.11);
\draw (5.44,-2.16)-- (5.4,-4.53);
\draw (3.34,-2.56)-- (3.52,-4.7);
\draw (3.52,-4.7)-- (5.4,-4.53);
\draw (5.4,-4.53)-- (7.06,-3.61);
\draw (5.4,-4.53)-- (4.26,-7.24);
\draw (3.52,-4.7)-- (4.26,-7.24);
\draw (7.19,-9.47)-- (4.26,-7.24);
\draw (4.26,-7.24)-- (2.03,-9.87);
\draw (7.19,-9.47)-- (2.03,-9.87);
\draw (2.03,-9.87)-- (-2.52,-7.55);
\draw (3.34,-2.56)-- (1.77,-4.05);
\draw (1.77,-4.05)-- (3.52,-4.7);
\draw (2.03,-9.87)-- (-0.38,-5.75);
\draw (1.77,-4.05)-- (-0.38,-5.75);
\draw (-0.33,-1.11)-- (0.41,-2.82);
\draw (0.41,-2.82)-- (1.77,-4.05);
\draw (-3,-1.64)-- (0.41,-2.82);
\draw (0.41,-2.82)-- (-0.38,-5.75);
\draw (-2.52,-7.55)-- (-0.38,-5.75);
\draw (-5.71,1.86)-- (-5.41,-3.3);
\draw (-5.41,-3.3)-- (-3,-1.64);
\draw (-5.41,-3.3)-- (-2.52,-7.55);
\draw (3.91,-0.18)-- (5.09,1.99);
\draw (3.91,-0.18)-- (6.41,0.07);
\draw (3.91,-0.18)-- (3.34,-2.56);
\draw (3.91,-0.18)-- (5.44,-2.16);
\draw (3.91,-0.18)-- (1.55,-0.37);
\draw (3.91,-0.18)-- (2.42,1.77);
\draw (0.89,9.3)-- (3.31,12.64);
\draw (6.01,9.61)-- (3.31,12.64);
\draw (13.71,-2.08)-- (15.21,-5.91);
\draw (15.21,-5.91)-- (11.4,-6.67);
\draw (-5.41,-3.3)-- (-6.66,-7.06);
\draw (-6.66,-7.06)-- (-2.52,-7.55);
\begin{scriptsize}
\fill [color=black] (0.89,9.3) circle (1.5pt);
\fill [color=black] (6.01,9.61) circle (1.5pt);
\fill [color=black] (1.33,6.41) circle (1.5pt);
\fill [color=black] (3.61,4.31) circle (1.5pt);
\fill [color=black] (6.19,6.59) circle (1.5pt);
\fill [color=black] (10.56,7.38) circle (1.5pt);
\fill [color=black] (-3.4,6.5) circle (1.5pt);
\fill [color=black] (1.86,3.7) circle (1.5pt);
\fill [color=black] (5.49,4.09) circle (1.5pt);
\fill [color=black] (10.04,3.57) circle (1.5pt);
\fill [color=black] (13.45,3) circle (1.5pt);
\fill [color=black] (7.11,3.13) circle (1.5pt);
\fill [color=black] (5.09,1.99) circle (1.5pt);
\fill [color=black] (2.42,1.77) circle (1.5pt);
\fill [color=black] (11,-1.11) circle (1.5pt);
\fill [color=black] (13.71,-2.08) circle (1.5pt);
\fill [color=black] (11.4,-6.67) circle (1.5pt);
\fill [color=black] (8.2,1.6) circle (1.5pt);
\fill [color=black] (8.55,-0.28) circle (1.5pt);
\fill [color=black] (6.41,0.07) circle (1.5pt);
\fill [color=black] (8.16,-2.03) circle (1.5pt);
\fill [color=black] (7.06,-3.61) circle (1.5pt);
\fill [color=black] (8.73,-5.45) circle (1.5pt);
\fill [color=black] (7.19,-9.47) circle (1.5pt);
\fill [color=black] (0.41,2.47) circle (1.5pt);
\fill [color=black] (-2.35,3.17) circle (1.5pt);
\fill [color=black] (5.44,-2.16) circle (1.5pt);
\fill [color=black] (-0.29,0.77) circle (1.5pt);
\fill [color=black] (1.55,-0.37) circle (1.5pt);
\fill [color=black] (3.34,-2.56) circle (1.5pt);
\fill [color=black] (-5.71,1.86) circle (1.5pt);
\fill [color=black] (-0.33,-1.11) circle (1.5pt);
\fill [color=black] (-3,-1.64) circle (1.5pt);
\fill [color=black] (5.4,-4.53) circle (1.5pt);
\fill [color=black] (3.52,-4.7) circle (1.5pt);
\fill [color=black] (4.26,-7.24) circle (1.5pt);
\fill [color=black] (2.03,-9.87) circle (1.5pt);
\fill [color=black] (-2.52,-7.55) circle (1.5pt);
\fill [color=black] (1.77,-4.05) circle (1.5pt);
\fill [color=black] (-0.38,-5.75) circle (1.5pt);
\fill [color=black] (0.41,-2.82) circle (1.5pt);
\fill [color=black] (-5.41,-3.3) circle (1.5pt);
\fill [color=black] (3.91,-0.18) circle (1.5pt);
\fill [color=ffffff] (3.31,12.64) circle (1.5pt);
\fill [color=ffffff] (15.21,-5.91) circle (1.5pt);
\fill [color=ffffff] (-6.66,-7.06) circle (1.5pt);
\draw (3.31,12.64) circle (1.5pt);
\draw (15.21,-5.91) circle (1.5pt);
\draw (-6.66,-7.06) circle (1.5pt);
\end{scriptsize}
\end{tikzpicture}
} 
\subfigure[]{\begin{tikzpicture}[line cap=round,line join=round,>=triangle 45,x=0.08cm,y=0.08cm]
\clip(-10.32,-14.16) rectangle (18.99,14.43);
\draw (0.89,9.3)-- (6.01,9.61);
\draw (0.89,9.3)-- (1.33,6.41);
\draw (1.33,6.41)-- (3.61,4.31);
\draw (3.61,4.31)-- (6.19,6.59);
\draw (6.01,9.61)-- (6.19,6.59);
\draw (6.19,6.59)-- (10.56,7.38);
\draw (10.56,7.38)-- (6.01,9.61);
\draw (0.89,9.3)-- (-3.4,6.5);
\draw (-3.4,6.5)-- (1.33,6.41);
\draw (1.33,6.41)-- (1.86,3.7);
\draw (1.86,3.7)-- (3.61,4.31);
\draw (3.61,4.31)-- (5.49,4.09);
\draw (6.19,6.59)-- (5.49,4.09);
\draw (10.56,7.38)-- (10.04,3.57);
\draw (10.56,7.38)-- (13.45,3);
\draw (10.04,3.57)-- (13.45,3);
\draw (10.04,3.57)-- (7.11,3.13);
\draw (5.49,4.09)-- (7.11,3.13);
\draw (5.49,4.09)-- (5.09,1.99);
\draw (5.09,1.99)-- (7.11,3.13);
\draw (2.42,1.77)-- (1.86,3.7);
\draw (2.42,1.77)-- (5.09,1.99);
\draw (13.45,3)-- (11,-1.11);
\draw (13.45,3)-- (13.71,-2.08);
\draw (11,-1.11)-- (13.71,-2.08);
\draw (13.71,-2.08)-- (11.4,-6.67);
\draw (7.11,3.13)-- (8.2,1.6);
\draw (10.04,3.57)-- (8.2,1.6);
\draw (8.2,1.6)-- (8.55,-0.28);
\draw (8.2,1.6)-- (6.41,0.07);
\draw (6.41,0.07)-- (8.55,-0.28);
\draw (8.55,-0.28)-- (11,-1.11);
\draw (8.55,-0.28)-- (8.16,-2.03);
\draw (11,-1.11)-- (8.16,-2.03);
\draw (8.16,-2.03)-- (7.06,-3.61);
\draw (8.16,-2.03)-- (8.73,-5.45);
\draw (7.06,-3.61)-- (8.73,-5.45);
\draw (11.4,-6.67)-- (8.73,-5.45);
\draw (8.73,-5.45)-- (7.19,-9.47);
\draw (11.4,-6.67)-- (7.19,-9.47);
\draw (1.86,3.7)-- (0.41,2.47);
\draw (0.41,2.47)-- (2.42,1.77);
\draw (-3.4,6.5)-- (-2.35,3.17);
\draw (-2.35,3.17)-- (0.41,2.47);
\draw (6.41,0.07)-- (5.44,-2.16);
\draw (5.44,-2.16)-- (7.06,-3.61);
\draw (-2.35,3.17)-- (-0.29,0.77);
\draw (0.41,2.47)-- (-0.29,0.77);
\draw (-0.29,0.77)-- (1.55,-0.37);
\draw (1.55,-0.37)-- (3.34,-2.56);
\draw (-3.4,6.5)-- (-5.71,1.86);
\draw (-2.35,3.17)-- (-5.71,1.86);
\draw (-0.29,0.77)-- (-0.33,-1.11);
\draw (-0.33,-1.11)-- (1.55,-0.37);
\draw (-5.71,1.86)-- (-3,-1.64);
\draw (-3,-1.64)-- (-0.33,-1.11);
\draw (5.44,-2.16)-- (5.4,-4.53);
\draw (3.34,-2.56)-- (3.52,-4.7);
\draw (3.52,-4.7)-- (5.4,-4.53);
\draw (5.4,-4.53)-- (7.06,-3.61);
\draw (5.4,-4.53)-- (4.26,-7.24);
\draw (3.52,-4.7)-- (4.26,-7.24);
\draw (7.19,-9.47)-- (4.26,-7.24);
\draw (4.26,-7.24)-- (2.03,-9.87);
\draw (7.19,-9.47)-- (2.03,-9.87);
\draw (2.03,-9.87)-- (-2.52,-7.55);
\draw (3.34,-2.56)-- (1.77,-4.05);
\draw (1.77,-4.05)-- (3.52,-4.7);
\draw (2.03,-9.87)-- (-0.38,-5.75);
\draw (1.77,-4.05)-- (-0.38,-5.75);
\draw (-0.33,-1.11)-- (0.41,-2.82);
\draw (0.41,-2.82)-- (1.77,-4.05);
\draw (-3,-1.64)-- (0.41,-2.82);
\draw (0.41,-2.82)-- (-0.38,-5.75);
\draw (-2.52,-7.55)-- (-0.38,-5.75);
\draw (-5.71,1.86)-- (-5.41,-3.3);
\draw (-5.41,-3.3)-- (-3,-1.64);
\draw (-5.41,-3.3)-- (-2.52,-7.55);
\draw (3.91,-0.18)-- (5.09,1.99);
\draw (3.91,-0.18)-- (6.41,0.07);
\draw (3.91,-0.18)-- (3.34,-2.56);
\draw (3.91,-0.18)-- (5.44,-2.16);
\draw (3.91,-0.18)-- (1.55,-0.37);
\draw (3.91,-0.18)-- (2.42,1.77);
\draw (0.89,9.3)-- (3.31,12.64);
\draw (6.01,9.61)-- (3.31,12.64);
\draw (13.71,-2.08)-- (15.21,-5.91);
\draw (15.21,-5.91)-- (11.4,-6.67);
\draw (-5.41,-3.3)-- (-6.66,-7.06);
\draw (-6.66,-7.06)-- (-2.52,-7.55);
\draw(3.94,-0.05) circle (12.71);
\begin{scriptsize}
\fill [color=black] (0.89,9.3) circle (1.5pt);
\fill [color=black] (6.01,9.61) circle (1.5pt);
\fill [color=black] (1.33,6.41) circle (1.5pt);
\fill [color=black] (3.61,4.31) circle (1.5pt);
\fill [color=black] (6.19,6.59) circle (1.5pt);
\fill [color=black] (10.56,7.38) circle (1.5pt);
\fill [color=black] (-3.4,6.5) circle (1.5pt);
\fill [color=black] (1.86,3.7) circle (1.5pt);
\fill [color=black] (5.49,4.09) circle (1.5pt);
\fill [color=black] (10.04,3.57) circle (1.5pt);
\fill [color=black] (13.45,3) circle (1.5pt);
\fill [color=black] (7.11,3.13) circle (1.5pt);
\fill [color=black] (5.09,1.99) circle (1.5pt);
\fill [color=black] (2.42,1.77) circle (1.5pt);
\fill [color=black] (11,-1.11) circle (1.5pt);
\fill [color=black] (13.71,-2.08) circle (1.5pt);
\fill [color=black] (11.4,-6.67) circle (1.5pt);
\fill [color=black] (8.2,1.6) circle (1.5pt);
\fill [color=black] (8.55,-0.28) circle (1.5pt);
\fill [color=black] (6.41,0.07) circle (1.5pt);
\fill [color=black] (8.16,-2.03) circle (1.5pt);
\fill [color=black] (7.06,-3.61) circle (1.5pt);
\fill [color=black] (8.73,-5.45) circle (1.5pt);
\fill [color=black] (7.19,-9.47) circle (1.5pt);
\fill [color=black] (0.41,2.47) circle (1.5pt);
\fill [color=black] (-2.35,3.17) circle (1.5pt);
\fill [color=black] (5.44,-2.16) circle (1.5pt);
\fill [color=black] (-0.29,0.77) circle (1.5pt);
\fill [color=black] (1.55,-0.37) circle (1.5pt);
\fill [color=black] (3.34,-2.56) circle (1.5pt);
\fill [color=black] (-5.71,1.86) circle (1.5pt);
\fill [color=black] (-0.33,-1.11) circle (1.5pt);
\fill [color=black] (-3,-1.64) circle (1.5pt);
\fill [color=black] (5.4,-4.53) circle (1.5pt);
\fill [color=black] (3.52,-4.7) circle (1.5pt);
\fill [color=black] (4.26,-7.24) circle (1.5pt);
\fill [color=black] (2.03,-9.87) circle (1.5pt);
\fill [color=black] (-2.52,-7.55) circle (1.5pt);
\fill [color=black] (1.77,-4.05) circle (1.5pt);
\fill [color=black] (-0.38,-5.75) circle (1.5pt);
\fill [color=black] (0.41,-2.82) circle (1.5pt);
\fill [color=black] (-5.41,-3.3) circle (1.5pt);
\fill [color=black] (3.91,-0.18) circle (1.5pt);
\fill [color=black] (3.31,12.64) circle (1.5pt);
\fill [color=black] (15.21,-5.91) circle (1.5pt);
\fill [color=black] (-6.66,-7.06) circle (1.5pt);
\end{scriptsize}
\end{tikzpicture}}
\subfigure[]{
\definecolor{black}{rgb}{0,0,0}
\begin{tikzpicture}[line cap=round,line join=round,>=triangle 45,x=0.075cm,y=0.075cm]
\clip(-18.51,-13.64) rectangle (12.78,15.3);
\draw (4.53,13.94)-- (11.8,2.26);
\draw (11.8,2.26)-- (5.19,-11.32);
\draw (4.53,13.94)-- (-10.44,13.94);
\draw (-10.44,13.94)-- (-17.64,2.12);
\draw (-17.64,2.12)-- (-10.22,-11.32);
\draw (5.19,-11.32)-- (-10.22,-11.32);
\draw (11.8,2.26)-- (4.53,13.94);
\draw (4.53,13.94)-- (11.8,2.26);
\draw (11.8,2.26)-- (7.18,-3.02);
\draw (5.19,-11.32)-- (7.18,-3.02);
\draw (4.53,13.94)-- (6.74,6.82);
\draw (6.74,6.82)-- (11.8,2.26);
\draw (6.74,6.82)-- (3.58,8.73);
\draw (4.53,13.94)-- (-2.96,12.4);
\draw (-2.96,12.4)-- (-10.44,13.94);
\draw (0.35,10.56)-- (-2.96,12.4);
\draw (-2.96,12.4)-- (-6.04,10.71);
\draw (0.35,10.56)-- (-6.04,10.71);
\draw (0.35,10.56)-- (3.58,8.73);
\draw (0.35,10.56)-- (0.79,8.14);
\draw (0.79,8.14)-- (3.58,8.73);
\draw (0.79,8.14)-- (0.35,5.35);
\draw (0.35,5.35)-- (-2.74,7.62);
\draw (0.79,8.14)-- (-2.74,7.62);
\draw (-2.74,7.62)-- (-6.7,8.06);
\draw (-6.7,8.06)-- (-6.04,10.71);
\draw (-6.7,8.06)-- (-9.42,8.73);
\draw (-6.04,10.71)-- (-9.42,8.73);
\draw (-9.42,8.73)-- (-12.5,6.82);
\draw (-12.5,6.82)-- (-10.44,13.94);
\draw (-12.5,6.82)-- (-17.64,2.12);
\draw (-17.64,2.12)-- (-12.43,-4.2);
\draw (-12.43,-4.2)-- (-10.22,-11.32);
\draw (-12.43,-4.2)-- (-12.5,-0.75);
\draw (-12.5,-0.75)-- (-8.9,-6.18);
\draw (-12.43,-4.2)-- (-8.9,-6.18);
\draw (-12.5,3.37)-- (-12.5,-0.75);
\draw (-12.5,6.82)-- (-12.5,3.37);
\draw (-9.42,8.73)-- (-12.5,3.37);
\draw (-12.5,-0.75)-- (-10.59,1.46);
\draw (-12.5,3.37)-- (-10.59,1.46);
\draw (-10.59,1.46)-- (-6.11,-2.43);
\draw (-6.11,-2.43)-- (-6.41,-5.45);
\draw (-8.9,-6.18)-- (-6.41,-5.45);
\draw (-5.67,-7.87)-- (-6.41,-5.45);
\draw (-8.9,-6.18)-- (-5.67,-7.87);
\draw (-5.67,-7.87)-- (-2.51,-9.48);
\draw (-5.67,-7.87)-- (0.35,-7.79);
\draw (-2.51,-9.48)-- (-10.22,-11.32);
\draw (-2.51,-9.48)-- (5.19,-11.32);
\draw (-2.51,-9.48)-- (0.35,-7.79);
\draw (0.35,-7.79)-- (4.17,-5.08);
\draw (4.17,-5.08)-- (1.52,-4.78);
\draw (0.35,-7.79)-- (1.52,-4.78);
\draw (7.18,-3.02)-- (4.17,-5.08);
\draw (1.52,-4.78)-- (-2.44,-4.93);
\draw (-2.44,-4.93)-- (-6.41,-5.45);
\draw (-2.44,-4.93)-- (-6.11,-2.43);
\draw (-2.44,-4.93)-- (1.01,-2.51);
\draw (1.01,-2.51)-- (-0.68,-1.26);
\draw (-0.68,-1.26)-- (3.65,-1.48);
\draw (3.65,-1.48)-- (1.01,-2.51);
\draw (-0.68,-1.26)-- (-4.42,-0.53);
\draw (-10.59,1.46)-- (-4.42,-0.53);
\draw (-4.42,-0.53)-- (-6.11,-2.43);
\draw (-4.42,-0.53)-- (-2.66,0.21);
\draw (-2.66,0.21)-- (-0.68,-1.26);
\draw (-2.66,0.21)-- (-7.21,2.48);
\draw (-2.66,0.21)-- (-2.66,1.6);
\draw (-2.66,1.6)-- (-7.21,2.48);
\draw (-7.21,2.48)-- (-12.5,3.37);
\draw (-7.21,2.48)-- (-7.14,4.03);
\draw (-7.14,4.03)-- (-12.5,3.37);
\draw (-7.14,4.03)-- (-5.97,5.5);
\draw (-4.2,4.17)-- (-7.14,4.03);
\draw (-5.97,5.5)-- (-4.2,4.17);
\draw (-5.97,5.5)-- (-6.7,8.06);
\draw (-5.97,5.5)-- (-2.74,7.62);
\draw (-4.2,4.17)-- (-0.68,4.25);
\draw (-0.68,4.25)-- (0.35,5.35);
\draw (-0.68,4.25)-- (4.83,1.9);
\draw (4.83,1.9)-- (0.35,5.35);
\draw (-0.68,4.25)-- (-2.59,3);
\draw (-4.2,4.17)-- (-2.59,3);
\draw (-2.59,3)-- (-2.66,1.6);
\draw (-2.59,3)-- (0.57,1.09);
\draw (0.57,1.09)-- (-2.66,1.6);
\draw (0.57,1.09)-- (3.65,-1.48);
\draw (1.01,-2.51)-- (1.52,-4.78);
\draw (3.65,-1.48)-- (7.18,0.06);
\draw (0.57,1.09)-- (7.18,0.06);
\draw (4.83,1.9)-- (7.18,0.06);
\draw (4.83,1.9)-- (7.1,3.22);
\draw (7.1,3.22)-- (7.18,0.06);
\draw (7.18,0.06)-- (4.17,-5.08);
\draw (7.18,0.06)-- (7.18,-3.02);
\draw (6.74,6.82)-- (7.1,3.22);
\draw (3.58,8.73)-- (7.1,3.22);
\begin{scriptsize}
\fill [color=black] (4.53,13.94) circle (1.5pt);
\fill [color=black] (11.8,2.26) circle (1.5pt);
\fill [color=black] (5.19,-11.32) circle (1.5pt);
\fill [color=black] (-10.44,13.94) circle (1.5pt);
\fill [color=black] (-17.64,2.12) circle (1.5pt);
\fill [color=black] (-10.22,-11.32) circle (1.5pt);
\fill [color=black] (7.18,-3.02) circle (1.5pt);
\fill [color=black] (6.74,6.82) circle (1.5pt);
\fill [color=black] (3.58,8.73) circle (1.5pt);
\fill [color=black] (-2.96,12.4) circle (1.5pt);
\fill [color=black] (0.35,10.56) circle (1.5pt);
\fill [color=black] (-6.04,10.71) circle (1.5pt);
\fill [color=black] (0.79,8.14) circle (1.5pt);
\fill [color=black] (0.35,5.35) circle (1.5pt);
\fill [color=black] (-2.74,7.62) circle (1.5pt);
\fill [color=black] (-6.7,8.06) circle (1.5pt);
\fill [color=black] (-9.42,8.73) circle (1.5pt);
\fill [color=black] (-12.5,6.82) circle (1.5pt);
\fill [color=black] (-12.43,-4.2) circle (1.5pt);
\fill [color=black] (-12.5,-0.75) circle (1.5pt);
\fill [color=black] (-8.9,-6.18) circle (1.5pt);
\fill [color=black] (-12.5,3.37) circle (1.5pt);
\fill [color=black] (-10.59,1.46) circle (1.5pt);
\fill [color=black] (-6.11,-2.43) circle (1.5pt);
\fill [color=black] (-6.41,-5.45) circle (1.5pt);
\fill [color=black] (-5.67,-7.87) circle (1.5pt);
\fill [color=black] (-2.51,-9.48) circle (1.5pt);
\fill [color=black] (0.35,-7.79) circle (1.5pt);
\fill [color=black] (4.17,-5.08) circle (1.5pt);
\fill [color=black] (1.52,-4.78) circle (1.5pt);
\fill [color=black] (-2.44,-4.93) circle (1.5pt);
\fill [color=black] (1.01,-2.51) circle (1.5pt);
\fill [color=black] (-0.68,-1.26) circle (1.5pt);
\fill [color=black] (3.65,-1.48) circle (1.5pt);
\fill [color=black] (-4.42,-0.53) circle (1.5pt);
\fill [color=black] (-2.66,0.21) circle (1.5pt);
\fill [color=black] (-7.21,2.48) circle (1.5pt);
\fill [color=black] (-2.66,1.6) circle (1.5pt);
\fill [color=black] (-7.14,4.03) circle (1.5pt);
\fill [color=black] (-5.97,5.5) circle (1.5pt);
\fill [color=black] (-4.2,4.17) circle (1.5pt);
\fill [color=black] (-0.68,4.25) circle (1.5pt);
\fill [color=black] (4.83,1.9) circle (1.5pt);
\fill [color=black] (-2.59,3) circle (1.5pt);
\fill [color=black] (0.57,1.09) circle (1.5pt);
\fill [color=black] (7.18,0.06) circle (1.5pt);
\fill [color=black] (7.1,3.22) circle (1.5pt);
\end{scriptsize}
\end{tikzpicture}
}

\caption{Five $C_4$-free planar graphs with $\delta = 4$.\label{fig20101011-1}}
\end{figure}

Figure \ref{fig20101011-1} illustrates five $C_4$-free planar graphs with minimum degree 4, for $n=30,36,44,46,47$ respectively (where the three white vertices of graph (c) are identified). Note that each planar graph in Figure \ref{fig20101011-1}(b),(d),(e) has at least one 6-face.

We begin to construct a new $C_4$-free planar graph with $n$ vertices and with minimum degree $\delta(G) = 4$ from one of the graphs illustrated in Figure \ref{fig20101011-1} (b), (d), (e).

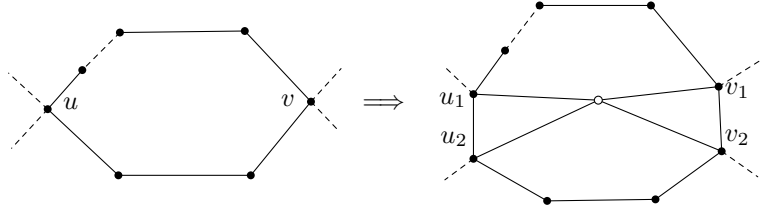
\begin{figure}[H]
\centering
\begin{tikzpicture}[line cap=round,line join=round,>=triangle 45,x=1.0cm,y=1.0cm]
\clip(-3.24,0.96) rectangle (7.5,4.06);
\draw (1.68,2.78) node[anchor=north west] {$\Longrightarrow$};
\draw (-1.42,3.54)-- (0.24,3.56);
\draw [dash pattern=on 2pt off 2pt] (-1.42,3.54)-- (-1.92,3.04);
\draw (-1.92,3.04)-- (-2.38,2.52);
\draw (-2.38,2.52)-- (-1.44,1.64);
\draw (-1.44,1.64)-- (0.32,1.64);
\draw (0.24,3.56)-- (1.12,2.62);
\draw (1.12,2.62)-- (0.32,1.64);
\draw (4.16,3.9)-- (5.64,3.9);
\draw (5.64,3.9)-- (6.54,2.82);
\draw (6.54,2.82)-- (6.58,1.96);
\draw (6.58,1.96)-- (5.7,1.32);
\draw (5.7,1.32)-- (4.26,1.3);
\draw (4.26,1.3)-- (3.28,1.86);
\draw (3.28,1.86)-- (3.28,2.72);
\draw (3.28,2.72)-- (3.7,3.3);
\draw [dash pattern=on 2pt off 2pt] (3.7,3.3)-- (4.16,3.9);
\draw (4.94,2.64)-- (3.28,2.72);
\draw (4.94,2.64)-- (3.28,1.86);
\draw (4.94,2.64)-- (6.54,2.82);
\draw (4.94,2.64)-- (6.58,1.96);
\draw [dash pattern=on 2pt off 2pt] (1.12,2.62)-- (1.52,3.06);
\draw (2.7,2.87) node[anchor=north west] {$u_1$};
\draw (2.72,2.35) node[anchor=north west] {$u_2$};
\draw (6.5,3) node[anchor=north west] {$v_1$};
\draw (6.5,2.37) node[anchor=north west] {$v_2$};
\draw (0.6,2.84) node[anchor=north west] {$v$};
\draw (-2.3,2.8) node[anchor=north west] {$u$};
\draw [dash pattern=on 2pt off 2pt] (1.12,2.62)-- (1.5,2.22);
\draw [dash pattern=on 2pt off 2pt] (-2.38,2.52)-- (-2.9,3);
\draw [dash pattern=on 2pt off 2pt] (-2.38,2.52)-- (-2.86,2);
\draw [dash pattern=on 2pt off 2pt] (3.28,2.72)-- (2.9,3.08);
\draw [dash pattern=on 2pt off 2pt] (3.28,1.86)-- (2.9,1.58);
\draw [dash pattern=on 2pt off 2pt] (6.54,2.82)-- (7.08,3.16);
\draw [dash pattern=on 2pt off 2pt] (6.58,1.96)-- (7.06,1.62);
\begin{scriptsize}
\fill [color=black] (-1.42,3.54) circle (1.5pt);
\fill [color=black] (0.24,3.56) circle (1.5pt);
\fill [color=black] (-1.92,3.04) circle (1.5pt);
\fill [color=black] (-2.38,2.52) circle (1.5pt);
\fill [color=black] (-1.44,1.64) circle (1.5pt);
\fill [color=black] (0.32,1.64) circle (1.5pt);
\fill [color=black] (1.12,2.62) circle (1.5pt);
\fill [color=black] (4.16,3.9) circle (1.5pt);
\fill [color=black] (5.64,3.9) circle (1.5pt);
\fill [color=black] (6.54,2.82) circle (1.5pt);
\fill [color=black] (6.58,1.96) circle (1.5pt);
\fill [color=black] (5.7,1.32) circle (1.5pt);
\fill [color=black] (4.26,1.3) circle (1.5pt);
\fill [color=black] (3.28,1.86) circle (1.5pt);
\fill [color=black] (3.28,2.72) circle (1.5pt);
\fill [color=black] (3.7,3.3) circle (1.5pt);
\fill [color=white] (4.94,2.64) circle (1.5pt);
\draw(4.94,2.64) circle (1.5pt);
\end{scriptsize}
\end{tikzpicture}
\caption{Operation (A)}\label{fig20101011}
\end{figure}

Take one $k$-face $f$ with $k\ge 6$, we construct a new planar graph $G^*$ by operation (A) which is illustrated in Figure \ref{fig20101011}. In this operation, we find two vertices $u, v$ which has distance 3 on $f$, then split $u$ and $v$ into two vertices $u_1, u_2$ and $v_1,v_2$ respectively. Finally we add a new vertex (the white vertex) inside $f$, and add edges from it to $u_1,u_2,v_1,v_2$ respectively. The resulting graph $G^*$ is $C_4$-free and with minimum degree $\delta = 4$ and with three more vertices. We can see that $G^*$ still has a $k$-face with $k\ge 6$ (one of the face incident with $u_1u_2$ or $v_1v_2$). So we can take operation (A) again on the $k$-face (with $k\ge 6$) on $G^*$, and therefore get a new $C_4$-free planar graph with  minimum degree $\delta = 4$ and with three more vertices than $G^{*}$.

In a result, if we take operation (A) recursively on graph (b) in Figure \ref{fig20101011-1}, we can construct $C_4$-free planar graph with $36+3t_1$ vertices and with minimum degree $\delta = 4$, where $t_1\ge 0$. Similarly,  if we take operation (A) recursively on graph (d) and (e)  in Figure \ref{fig20101011-1} respectively, we can get graphs with $46+3t_2$ and $47+3t_3$ vertices, where $t_2,t_3\ge 0$.  This complete the proof of (i) in Theorem \ref{thm20101011}. 

(ii) Let $G$ be a $C_4$-free planar graph of order $n$ and  with $\delta(G)=\delta(n,C_4)$. If $31\leq n\leq 33$, then by Lemma \ref{lem130420}, we have $\varepsilon(G)\le 2n-1$, this implies that  $\delta(n,C_4)\leq 3$.  If $n=34,35,37,38$, then by Lemma \ref{lem130420} again, we have $M(n,C_4)=2n$ for $n=34,35,37,38$, this implies that $\delta(G)\leq 4$. If $\delta(G)=4$, then $G$ is 4-regular and $\varepsilon(G)=M(n,C_4)=2n$. By Theorem \ref{thm3}, $\varepsilon(G)=M(34,C_4)=68$ if and only if $\tau(G)=2$ and $f_6=\cdots=f_r=0$; $\varepsilon(G)=M(35,C_4)=70$ if and only if $\tau(G)=1$, $f_6=1$ and $f_7=\cdots=f_r=0$;
$\varepsilon(G)=M(37,C_4)=74$ if and only if $\tau(G)=2$, $f_6=1$ and $f_7=\cdots=f_r=0$; $\varepsilon(G)=M(38,C_4)=76$ if and only if $\tau(G)=1$, $f_7=1$ and $f_6=f_8=\cdots=f_r=0$, or $\tau(G)=4$ and $f_6=\cdots=f_r=0$, or $\tau(G)=1$, $f_6=2$ and $f_7=\cdots=f_r=0$. But each of the above cases contradicts the facts of Lemma \ref{lem20101012}.

If $n=40$,  suppose on the contrary that there is a $C_4$-free planar graph of order $n=40$ with $\delta(G)\ge 4$.  By Theorem \ref{thm3}, we get that $40\le\varepsilon(G)\le 81$. 

Suppose first that $\varepsilon(G)=80$, then $G$ is 4-regular, and further more, there are only three possibilities to consider: (a) $\tau(G)=2,f_7=1$ and $f_6=f_8=\cdots=f_r=0$; (b) $\tau(G)=2, f_6=2$ and $f_7=f_8=\cdots=f_r=0$; (c) $\tau(G)=5$ and $f_6=f_7=\cdots=f_r=0$. By Lemma \ref{lem20101012}, the first two possibilities can not happen. 

So we assume that (c) holds. 
By Euler's formula, $G$ has 17 pentagons and 25 triangles. Since $\tau(G)=5$ and by Lemma \ref{lem20101012}, $\Gamma(G)$ induces a 5-cycle $C$ in $G$.  By Lemma \ref{lem20111029-1},  $C$ is a 5-face. Consider the vertex-edge-dual $G^*$ of $G$, since $G$ is $C_4$-free,  $G^*$ is a triangulation of 17 vertices with degree sequence $5^{12}6^5$; and furthermore, since $\Gamma(G)$ induces a 5-face,  $G^*$ has the property that there is a vertex of degree 5 which is adjacent to every vertex of degree 6.  By checking the graphs in Figure \ref{fig20111026-1}, none of them has the above property, a contradiction.
 So in the following, we suppose that $\varepsilon(G)=81$. 
 
 By Theorem \ref{thm3}, this is possible only if $\tau(G)=0, f_6=1$, and $f_7=f_8=\cdots =0$.
In this case $G$ has 15 pentagons, 27 triangles and one 6-face. By the definition of vertex-edge-dual $G^*$ of $G$, we see that $G^*$ is a triangulation on 16 vertices with $\delta(G^*)\ge 5$. Since $\delta(G)=4, \tau(G)=0$ and $\varepsilon(G)=81$, the degree sequence of $G$ is exactly $4^{39}6^1$. Let $f$ be the 6-face and $v$ be the vertex of degree 6 in $G$, let $f_1$, $f_2$, $f_3$ be the nontriangle faces which are incident to $v$. If $f$ is not adjacent to $v$, then by the rule of the construction of $G^*$, $G^*$ have exactly four vertices of degree 6 (corresponding to $f,f_1,f_2,f_3$) and three of which (corresponding to $f_1, f_2,f_3$)  form a triangle in $G^*$; If $f$ is adjacent to $v$, then $G^*$ have exactly two vertices of degree 6 and two vertices of degree 7. By Fact \ref{lem20111024}, there are all together 3 non-isomorphic triangulations on 16 vertices with minimum degree 5, and none of them has the above property, a contradiction. 

Therefore, we have the conclusion that $\delta(40,C_4)\le 3$. 

If $n=41$, suppose on the contrary that there is a $C_4$-free planar graph of order $n=41$ with $\delta(G)\ge 4$. By Theorem \ref{thm3}, we have that $82\le\varepsilon(G)\le 83$.  If $\varepsilon(G)=82$, then $G$ is 4-regular, and by Theorem \ref{thm3}, there are four possibilities to consider:

(1) $\tau(G)=1, f_6=f_7=1, f_8=f_9=\cdots=0$;

(2) $\tau(G)=1, f_6=3, f_7=f_8=\cdots=0$;

(3) $\tau(G)=1, f_8=1, f_6=f_7=f_9=\cdots=0$;

(4) $\tau(G)=4, f_6=1, f_7=f_8=\cdots=0$.

\begin{figure}[H]
\centering
\definecolor{ffffff}{rgb}{1,1,1}
\begin{tikzpicture}[line cap=round,line join=round,>=triangle 45,x=1.0cm,y=1.0cm]
\clip(-1.68,-0.7) rectangle (3.86,4.3);
\draw (0.92,2.94)-- (0.92,1.46);
\draw (0.92,2.94)-- (2.22,2.94);
\draw (2.22,2.94)-- (2.96,2.16);
\draw (2.96,2.16)-- (2.22,1.46);
\draw (0.92,1.46)-- (2.22,1.46);
\draw (0.92,2.94)-- (-0.48,2.94);
\draw (-0.48,2.94)-- (-1.14,2.18);
\draw (-1.14,2.18)-- (-0.44,1.46);
\draw (0.92,1.46)-- (-0.44,1.46);
\draw (0.92,2.94)-- (0.56,3.76);
\draw (0.92,2.94)-- (1.26,3.76);
\draw (2.22,1.46)-- (3.1,1.44);
\draw (2.96,2.16)-- (3.1,1.44);
\draw (0.92,1.46)-- (0.92,0.36);
\draw (0.92,0.36)-- (1.64,-0.14);
\draw (1.64,-0.14)-- (2.22,0.28);
\draw (2.22,1.46)-- (2.22,0.28);
\draw (2.22,1.46)-- (2.86,0.54);
\draw (2.86,0.54)-- (2.22,0.28);
\draw (2.22,2.94)-- (2.56,3.2);
\draw (2.96,2.16)-- (3.08,2.52);
\draw (-0.48,2.94)-- (-0.78,3.48);
\draw (-0.48,2.94)-- (-1.02,3.16);
\draw (-1.14,2.18)-- (-1.52,2.34);
\draw (-1.14,2.18)-- (-1.5,1.96);
\draw (-0.44,1.46)-- (-0.8,1.18);
\draw (0.92,0.36)-- (0.68,-0.12);
\draw (1.64,-0.14)-- (1.42,-0.54);
\draw (1.64,-0.14)-- (1.82,-0.56);
\draw (2.22,0.28)-- (2.38,-0.12);
\draw (2.86,0.54)-- (3.28,0.18);
\draw (2.86,0.54)-- (3.42,0.48);
\draw (3.1,1.44)-- (3.66,1.4);
\draw (3.1,1.44)-- (3.48,1.04);
\draw (1.26,3.76)-- (1.5,4.12);
\draw (0.56,3.76)-- (0.5,4.16);
\draw (0.56,3.76)-- (0.24,4);
\draw (-0.44,1.46)-- (0.92,0.36);
\draw (1.26,3.76)-- (2.22,2.94);
\draw (0.56,3.76)-- (-0.48,2.94);
\draw (1.26,3.76)-- (1.14,4.14);
\draw (0.9,3.04) node[anchor=north west] {$u$};
\draw (2,1.98) node[anchor=north west] {$v$};
\draw (0.88,1.98) node[anchor=north west] {$w$};
\draw (0.46,2.42) node[anchor=north west] {$e$};
\draw (1.34,1.5) node[anchor=north west] {$f$};
\begin{scriptsize}
\fill [color=ffffff] (0.92,2.94) circle (1.5pt);
\draw (0.92,2.94) circle (1.5pt);
\fill [color=black] (0.92,1.46) circle (1.5pt);
\fill [color=black] (2.22,2.94) circle (1.5pt);
\fill [color=black] (2.96,2.16) circle (1.5pt);
\fill [color=ffffff] (2.22,1.46) circle (1.5pt);
\draw (2.22,1.46) circle (1.5pt);
\fill [color=black] (-0.48,2.94) circle (1.5pt);
\fill [color=black] (-1.14,2.18) circle (1.5pt);
\fill [color=black] (-0.44,1.46) circle (1.5pt);
\fill [color=black] (0.56,3.76) circle (1.5pt);
\fill [color=black] (1.26,3.76) circle (1.5pt);
\fill [color=black] (3.1,1.44) circle (1.5pt);
\fill [color=black] (0.92,0.36) circle (1.5pt);
\fill [color=black] (1.64,-0.14) circle (1.5pt);
\fill [color=black] (2.22,0.28) circle (1.5pt);
\fill [color=black] (2.86,0.54) circle (1.5pt);
\end{scriptsize}
\end{tikzpicture}
\caption{Local structure of $G$.\label{fig20111028}}
\end{figure}
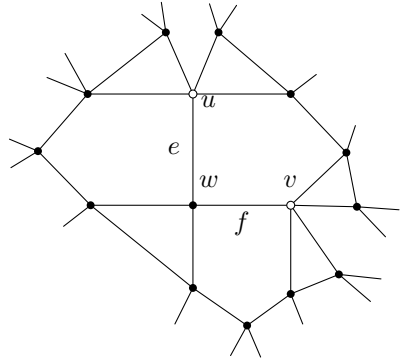

But all these cases contradicts Lemma \ref{lem20101012}. So we have that $\varepsilon(G)=83$. By Theorem \ref{thm3}, this can happen only if $\tau(G)=2, f_6=f_7=\cdots=0$.  Furthermore, since $\delta(G)\ge 4$, the degree sequence of $G$ is either $4^{40}6^1$ or $4^{39}5^2$. Assume first that the degree sequence of $G$ is 
$4^{40}6^1$. If there is an edge $e\in \Gamma(G)$  (say $e=uv$)  such that $d_G(u)=d_G(v)=4$,   by Lemma \ref{lem20101012} (iii),  we have that $\tau(G)\ge 3$, which contradicts that $\tau(G)=2$; If the two edges $e,f \in \Gamma(G)$ satisfy that $e=wu, f=wv$ with $d_G(u)=d_G(v)=4$ and $d_G(w)=6$, by Lemma \ref{lem20101012} (iii) again,  we have that $\tau(G)\ge 4$, which is a contradiction too. So we assume that the degree sequence of $G$ is $4^{39}5^2$.  Let $u,v$ in $G$ such that $d_G(u)=d_G(v)=5$ and let $\Gamma(G)=\{e,f\}$.  By Lemma \ref{lem20101012}, it suffices to consider the case that $e,f$ have a vertex in common and are incident with $u,v$ respectively (Figure \ref{fig20111028}). 

By Euler's Formula, $G$ has 27 triangles and 17 pentagons. Consider the vertex-edge-dual $G^*$ of $G$. It is a triangulation on 17 vertices with degree sequence $5^{14}6^3$, this is impossible since any triangulation on 17 vertices has $3\times 17-6=45$ edges. 

Therefore, we have the conclusion that $\delta(41,C_4)\le 3$. 

If $n=43$, suppose on the contrary that there is a $C_4$-free planar graph of order $n=43$ with $\delta(G)=4$.  By Theorem \ref{thm3}, we have $86\le\varepsilon(G)\le 87$.

If $\varepsilon(G)=86$, then $G$ is 4-regular since $\delta(G)=4$. By Theorem \ref{thm3}, there are four possibilities to consider:

(1) $\tau(G)=2, f_8=1, f_6=f_7=f_9=\cdots=0$;

(2) $\tau(G)=2, f_6=3, f_7=f_8=\cdots=0$;

(3) $\tau(G)=2, f_6=f_7=1, f_8=f_9=\cdots=0$;

(4) $\tau(G)=5, f_6=1, f_7=f_8\cdots=0$;

By Lemma \ref{lem20101012}, the first three cases can not happen. So we assume that $\tau(G)=5, f_6=1, f_7=f_8\cdots=0$. In this case, $G$ is 4-regular, with 27 triangles, 17 pentagons and one hexagon. 

Let $\Gamma(G)=\{e_1,e_2,e_3,e_4,e_5\}$. By Lemma \ref{lem20101012} and \ref{lem20111029-1}, $\Gamma(G)$ induces a  5-face in $G$. Consider the vertex-edge-dual  $G^*$ of $G$,  it is a triangulation on 18 vertices with degree sequence $5^{12}6^6$ or $5^{13}6^47^1$. Furthermore, $G^*$ has the additional property: let $T$ be the vertex set consisting of all vertices with degrees at least 6 in $G^*$, then there is a five cycle in the subgraph induced by $T$ in $G^*$. By checking the graphs in Fact \ref{lem20111029}, none of them has that property, a contradiction.

Now we assume $\varepsilon(G)=87$. By Theorem \ref{thm3}, we shall only consider the following three cases:

(1) $\tau(G)=0, f_6=2, f_7=f_8\cdots=0$;

(2) $\tau(G)=0, f_7=1, f_6=f_8=f_9=\cdots=0$;

(3) $\tau(G)=3, f_6=f_7=\cdots=0$.

For case (1), by Euler's formula, we have $f_5=15, f_6=2$.  Since $\delta(G)=4$ and $\varepsilon(G)=87$, the degree sequence of $G$ is $4^{42}6^1$ or $4^{41}5^2$. Since $G$ is $C_4$-free and $\tau(G)=0$, there is no vertex of degree 5 in $G$, so the degree sequence of $G$ is $4^{42}6^1$. Consider the vertex-edge-dual $G^*$ of $G$, it is a triangulation on 17 vertices; Furthermore, it has at least 3 vertices of degree at least 6, and three of them form a triangle in $G^*$.  By checking the graphs in Fact \ref{lem20111026}, we see that none of which has the above property, a contradiction.

For case (2), a similar argument as in case (i) will deduce a contradiction.

For case (3), by Euler's formula, we have $f_5=18$.  Since $\delta(G)=4$ and $\varepsilon(G)=87$, the degree sequence of $G$ is $4^{42}6^1$ or $4^{41}5^2$. 
Since $\tau(G)=3$, let $\Gamma(G)=\{e_1,e_2,e_3\}$. We shall consider the following four cases:

{\bf Case 1.} $\Gamma(G)$ forms a matching in $G$.

Since there are at most two vertices of degrees at least 5, there is an edge in $\Gamma(G)$ (say $e_1$) so that the degrees of both endpoints of $e_1$ are four. By Lemma \ref{lem20101012}, we have $\tau(G)\ge 5$, which contradicts that $\tau(G)=3$.

{\bf Case 2.} $\Gamma(G)$ induces two disjoint paths. 

If the degree sequence of $G$ is $4^{42}6^1$, then there must be an edge in $\Gamma(G)$ such that each endpoint of which has degree 4 in $G$, this implies by  Lemma \ref{lem20101012} that $\tau(G)\ge 4$, which contradicts that $\tau(G)=3$. 

If the degree sequence of $G$ is $4^{41}5^2$, let $P_1=u_1u_2u_3$, $P_2=v_1v_2$ be the two disjoint paths induced by $\Gamma(G)$ respectively, and let $u,v$ be the two vertices of degrees 5 in $G$. If $d_G(u_2)=4$, then at least two vertices of $\{u_1,u_3,v_1,v_2\}$ have degrees 4 in $G$, this implies by  Lemma \ref{lem20101012} that $\tau(G)\ge 5$, which contradicts that $\tau(G)=3$. If $d_G(u_2)=5$, then at east one vertex of $v_1,v_2$ has degree 4 in $G$, this implies by  Lemma \ref{lem20101012} that $\tau(G)\ge 4$, which contradicts again that $\tau(G)=3$.

{\bf Case 3.} $\Gamma(G)$ induces a path of length 3 in $G$.

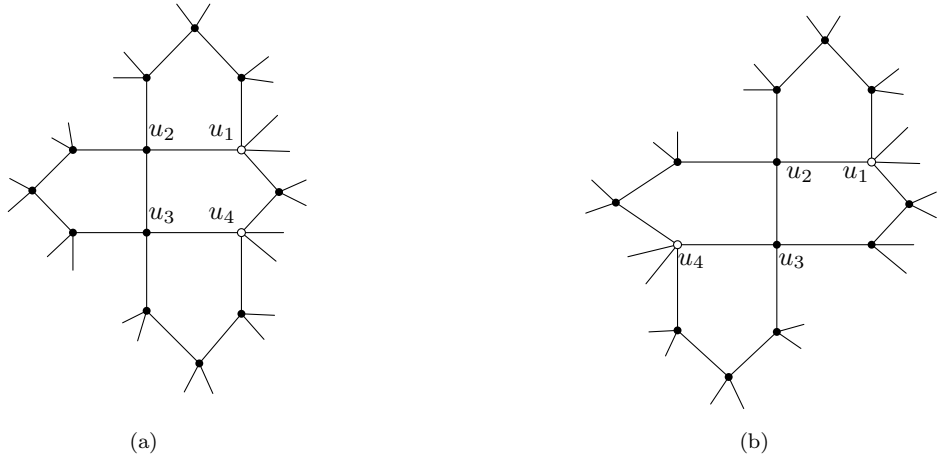
\begin{figure}[H]
\centering
\subfigure[]{\definecolor{ffffff}{rgb}{1,1,1}
\begin{tikzpicture}[line cap=round,line join=round,>=triangle 45,x=1.0cm,y=1.0cm]
\clip(-1.16,-0.7) rectangle (3.12,4.94);
\draw (2.16,2.84)-- (0.9,2.84);
\draw (0.9,2.84)-- (0.9,1.74);
\draw (0.9,1.74)-- (2.16,1.74);
\draw (2.66,2.28)-- (2.16,1.74);
\draw (2.66,2.28)-- (2.16,2.84);
\draw (0.9,2.84)-- (0.9,3.8);
\draw (2.16,2.84)-- (2.16,3.8);
\draw (0.9,3.8)-- (1.54,4.46);
\draw (1.54,4.46)-- (2.16,3.8);
\draw (2.16,1.74)-- (2.16,0.66);
\draw (0.9,1.74)-- (0.9,0.7);
\draw (2.16,0.66)-- (1.6,0);
\draw (0.9,0.7)-- (1.6,0);
\draw (0.9,2.84)-- (-0.08,2.84);
\draw (-0.08,2.84)-- (-0.62,2.3);
\draw (-0.62,2.3)-- (-0.08,1.74);
\draw (-0.08,1.74)-- (0.9,1.74);
\draw (2.16,2.84)-- (2.64,3.3);
\draw (2.16,2.84)-- (2.8,2.82);
\draw (2.16,1.74)-- (2.62,1.36);
\draw (2.16,1.74)-- (2.72,1.74);
\draw (1.6,3.3) node[anchor=north west] {$u_1$};
\draw (0.8,3.3) node[anchor=north west] {$u_2$};
\draw (0.8,2.2) node[anchor=north west] {$u_3$};
\draw (1.6,2.2) node[anchor=north west] {$u_4$};
\draw (2.16,3.8)-- (2.52,4.08);
\draw (2.16,3.8)-- (2.58,3.74);
\draw (1.54,4.46)-- (1.3,4.78);
\draw (1.54,4.46)-- (1.74,4.78);
\draw (0.9,3.8)-- (0.6,4.14);
\draw (0.9,3.8)-- (0.46,3.8);
\draw (-0.08,2.84)-- (-0.36,3);
\draw (-0.08,2.84)-- (-0.14,3.2);
\draw (-0.62,2.3)-- (-0.92,2.46);
\draw (-0.62,2.3)-- (-0.94,2.02);
\draw (-0.08,1.74)-- (-0.42,1.42);
\draw (-0.08,1.74)-- (-0.08,1.24);
\draw (0.9,0.7)-- (0.58,0.54);
\draw (0.9,0.7)-- (0.78,0.3);
\draw (1.6,0)-- (1.4,-0.38);
\draw (1.6,0)-- (1.78,-0.4);
\draw (2.16,0.66)-- (2.46,0.32);
\draw (2.16,0.66)-- (2.6,0.64);
\draw (2.66,2.28)-- (3.02,2.44);
\draw (2.66,2.28)-- (3.02,2.1);
\begin{scriptsize}
\fill [color=ffffff] (2.16,2.84) circle (1.5pt);
\draw (2.16,2.84) circle (1.5pt);
\fill [color=black] (0.9,2.84) circle (1.5pt);
\fill [color=black] (0.9,1.74) circle (1.5pt);
\fill [color=ffffff] (2.16,1.74) circle (1.5pt);
\draw (2.16,1.74) circle (1.5pt);
\fill [color=black] (2.66,2.28) circle (1.5pt);
\fill [color=black] (0.9,3.8) circle (1.5pt);
\fill [color=black] (2.16,3.8) circle (1.5pt);
\fill [color=black] (1.54,4.46) circle (1.5pt);
\fill [color=black] (2.16,0.66) circle (1.5pt);
\fill [color=black] (0.9,0.7) circle (1.5pt);
\fill [color=black] (1.6,0) circle (1.5pt);
\fill [color=black] (-0.08,2.84) circle (1.5pt);
\fill [color=black] (-0.62,2.3) circle (1.5pt);
\fill [color=black] (-0.08,1.74) circle (1.5pt);
\end{scriptsize}
\end{tikzpicture}}
\hskip 1.3in
\subfigure[]{\definecolor{ffffff}{rgb}{1,1,1}
\begin{tikzpicture}[line cap=round,line join=round,>=triangle 45,x=1.0cm,y=1.0cm]
\clip(-1.72,-0.54) rectangle (3.16,4.96);
\draw (2.16,2.84)-- (0.9,2.84);
\draw (0.9,2.84)-- (0.9,1.74);
\draw (0.9,1.74)-- (2.16,1.74);
\draw (2.66,2.28)-- (2.16,1.74);
\draw (2.66,2.28)-- (2.16,2.84);
\draw (0.9,2.84)-- (0.9,3.8);
\draw (2.16,2.84)-- (2.16,3.8);
\draw (0.9,3.8)-- (1.54,4.46);
\draw (1.54,4.46)-- (2.16,3.8);
\draw (-0.42,1.74)-- (0.9,1.74);
\draw (2.16,2.84)-- (2.64,3.3);
\draw (2.16,2.84)-- (2.8,2.82);
\draw (2.16,1.74)-- (2.62,1.36);
\draw (2.16,1.74)-- (2.72,1.74);
\draw (1.65,2.9) node[anchor=north west] {$u_1$};
\draw (0.9,2.9) node[anchor=north west] {$u_2$};
\draw (0.8,1.75) node[anchor=north west] {$u_3$};
\draw (-0.54,1.75) node[anchor=north west] {$u_4$};
\draw (2.16,3.8)-- (2.52,4.08);
\draw (2.16,3.8)-- (2.58,3.74);
\draw (1.54,4.46)-- (1.3,4.78);
\draw (1.54,4.46)-- (1.74,4.78);
\draw (0.9,3.8)-- (0.6,4.14);
\draw (0.9,3.8)-- (0.46,3.8);
\draw (2.66,2.28)-- (3.02,2.44);
\draw (2.66,2.28)-- (3.02,2.1);
\draw (-0.42,1.74)-- (-0.42,0.6);
\draw (0.9,1.74)-- (0.9,0.58);
\draw (-0.42,0.6)-- (0.26,-0.02);
\draw (0.9,0.58)-- (0.26,-0.02);
\draw (0.9,2.84)-- (-0.42,2.84);
\draw (-0.42,2.84)-- (-1.24,2.3);
\draw (-1.24,2.3)-- (-0.42,1.74);
\draw (-0.42,1.74)-- (-1.08,1.58);
\draw (-0.42,1.74)-- (-0.84,1.22);
\draw (-0.42,2.84)-- (-0.42,3.24);
\draw (-0.42,2.84)-- (-0.74,3.1);
\draw (-1.24,2.3)-- (-1.56,2.6);
\draw (-1.24,2.3)-- (-1.64,2.16);
\draw (0.9,0.58)-- (1.22,0.36);
\draw (0.9,0.58)-- (1.26,0.68);
\draw (0.26,-0.02)-- (0.46,-0.44);
\draw (0.26,-0.02)-- (0.02,-0.38);
\draw (-0.42,0.6)-- (-0.8,0.58);
\draw (-0.42,0.6)-- (-0.58,0.26);
\begin{scriptsize}
\fill [color=ffffff] (2.16,2.84) circle (1.5pt);
\draw (2.16,2.84) circle (1.5pt);
\fill [color=black] (0.9,2.84) circle (1.5pt);
\fill [color=black] (0.9,1.74) circle (1.5pt);
\fill [color=black] (2.16,1.74) circle (1.5pt);
\fill [color=black] (2.66,2.28) circle (1.5pt);
\fill [color=black] (0.9,3.8) circle (1.5pt);
\fill [color=black] (2.16,3.8) circle (1.5pt);
\fill [color=black] (1.54,4.46) circle (1.5pt);
\fill [color=ffffff] (-0.42,1.74) circle (1.5pt);
\draw (-0.42,1.74) circle (1.5pt);
\fill [color=black] (-0.42,0.6) circle (1.5pt);
\fill [color=black] (0.9,0.58) circle (1.5pt);
\fill [color=black] (0.26,-0.02) circle (1.5pt);
\fill [color=black] (-0.42,2.84) circle (1.5pt);
\fill [color=black] (-1.24,2.3) circle (1.5pt);
\end{scriptsize}
\end{tikzpicture}}
\caption{Local structure of $G$.\label{fig20111030}}
\end{figure}

Let $P=u_1u_2u_3u_4$ be the path induced by $\Gamma(G)$ in $G$. If $d_G(u_1)=4$ or $d_G(u_4)=4$, then $\tau(G)\ge 4$ by  Lemma \ref{lem20101012}, which contradicts that $\tau(G)=3$. This implies that the degree sequence of $G$ is exactly $4^{41}5^2$, and that $d_G(u_1)=d_G(u_4)=5$.  Then $G$ must have one of the following structure (Figure \ref{fig20111030} (a) or (b)).  Now we consider the vertex-edge-dual $G^*$ of $G$, note that $G$ has exactly 18 pentagons and no more faces of length at least 6,  we can see that in both cases $G^*$ are triangulations on 18 vertices with degree sequence $5^{14}6^4$, but this impossible since $G^{*}$ is a triangulation on 18 vertices. 

{\bf Case 4.} $\Gamma(G)$ induces a 3-cycle in $G$.

In this case $\Gamma(G)$ must induces a separating triangle in $G$. By a similar argument as in the proof of Lemma \ref{lem20111029-1} will deduce a contradiction.

Therefore, we have the conclusion that for each $n\in B=\{31, 32, 33, 34, 35, 37, 38, 40, 41, 43\}$, $\delta(n,C_4)\le 3$. 

\begin{figure}[H]
\centering
\begin{tikzpicture}[line cap=round,line join=round,>=triangle 45,x=1.0cm,y=1.0cm]
\clip(-2.3,2.64) rectangle (4.24,5.24);
\draw (-2,5)-- (-1,4);
\draw (-1,4)-- (0,3);
\draw (-1,4)-- (-2,3);
\draw (-1,4)-- (0,5);
\draw (2.5,4)-- (3.5,3.98);
\draw (2.5,4)-- (2,5);
\draw (2.5,4)-- (2,3);
\draw (3.5,3.98)-- (4,5);
\draw (3.5,3.98)-- (4,3);
\draw (-1.14,5.2) node[anchor=north west] {$f$};
\draw (-1.1,3.4) node[anchor=north west] {$g$};
\draw (2.9,5.14) node[anchor=north west] {$f$};
\draw (2.88,3.44) node[anchor=north west] {$g$};
\draw (-0.86,4.2) node[anchor=north west] {$v$};
\draw (1.88,4.2) node[anchor=north west] {$v_1$};
\draw (3.48,4.2) node[anchor=north west] {$v_2$};
\draw (0.74,4.2) node[anchor=north west] {$\Longrightarrow$};
\begin{scriptsize}
\fill [color=black] (-2,5) circle (1.5pt);
\fill [color=black] (-1,4) circle (1.5pt);
\fill [color=black] (0,3) circle (1.5pt);
\fill [color=black] (-2,3) circle (1.5pt);
\fill [color=black] (0,5) circle (1.5pt);
\fill [color=black] (2.5,4) circle (1.5pt);
\fill [color=black] (3.5,3.98) circle (1.5pt);
\fill [color=black] (2,5) circle (1.5pt);
\fill [color=black] (2,3) circle (1.5pt);
\fill [color=black] (4,5) circle (1.5pt);
\fill [color=black] (4,3) circle (1.5pt);
\end{scriptsize}
\end{tikzpicture}
\caption{Operation (B).\label{operation(B)}}
\end{figure}
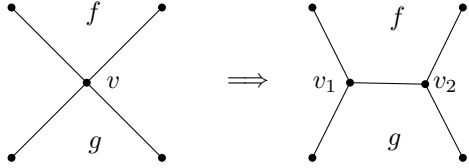

Now it remains to show that for each $n\in B=\{31, 32, 33, 34, 35, 37, 38, 40, 41, 43\}$, there is a $C_4$-free planar graph $G$ with $\delta(G)=3$.  We begin with the graph shown in Figure \ref{fig20101011-1} (a), it is a $C_4$-free 4-regular planar graph on 30 vertices. Each time we take one vertex $v$ with degree 4 in $G$, and make Operation (B) as shown in Figure \ref{operation(B)} (where $f$ and $g$ are faces of length at least 5). In this operation, the vertex $v$ is split to two vertices $v_1$ and $v_2$, then add an edge between $v_1$ and $v_2$. In this way, we get a new $C_4$-free planar graph with $\delta=3$ with one more vertex. If we make the Operation (B) $n-30$ times (note that $n\le 43$, each time we can always find a vertex of degree 4 in the new graph),  we finally get a $C_4$-free planar graph on $n$ vertices with $\delta=3$. 

\begin{figure}[H]
\centering
\definecolor{qqqqff}{rgb}{0,0,0}
\begin{tikzpicture}[line cap=round,line join=round,>=triangle 45,x=0.5cm,y=0.5cm]
\clip(-2.15,0.08) rectangle (5.18,5.64);
\draw (-0.44,2.76)-- (1.52,2.76);
\draw (-0.44,2.76)-- (-1.7,0.36);
\draw (1.54,0.36)-- (-1.7,0.36);
\draw (4.8,5.2)-- (4.8,0.36);
\draw (4.8,5.2)-- (3.48,4.1);
\draw (3.48,4.1)-- (3.48,2.78);
\draw (1.52,2.76)-- (3.48,2.78);
\draw (1.54,0.36)-- (4.8,0.36);
\draw (-1.7,5.2)-- (-0.44,2.76);
\draw (-1.7,5.2)-- (3.48,4.1);
\draw (-1.7,5.2)-- (4.8,5.2);
\draw (-1.7,5.2)-- (-1.7,0.36);
\draw (3.48,1.36)-- (1.52,2.76);
\draw (3.48,2.78)-- (3.48,1.36);
\draw (3.48,1.36)-- (1.54,0.36);
\draw (3.48,1.36)-- (4.8,0.36);
\begin{scriptsize}
\fill [color=qqqqff] (-0.44,2.76) circle (1.5pt);
\fill [color=qqqqff] (-1.7,0.36) circle (1.5pt);
\fill [color=qqqqff] (4.8,5.2) circle (1.5pt);
\fill [color=qqqqff] (3.48,4.1) circle (1.5pt);
\fill [color=qqqqff] (1.52,2.76) circle (1.5pt);
\fill [color=qqqqff] (1.54,0.36) circle (1.5pt);
\fill [color=qqqqff] (4.8,0.36) circle (1.5pt);
\fill [color=qqqqff] (3.48,2.78) circle (1.5pt);
\fill [color=qqqqff] (-1.7,5.2) circle (1.5pt);
\fill [color=qqqqff] (3.48,1.36) circle (1.5pt);
\end{scriptsize}
\end{tikzpicture}
\caption{A $C_4$-free planar graph on 10 vertices with $\delta =3$.\label{C4free10-m3}}
\end{figure}
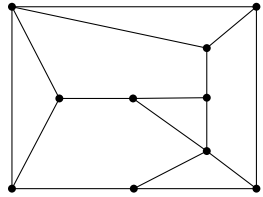

If $10\le n\le 29$, by corollary \ref{coro20101010}, $\delta(n,C_4)\le 3$. So it suffices to show the existence of $C_4$-free planar graph on $n$ vertices with minimum degree 3. We begin with the $C_4$-free planar graph $G$ on 10 vertices with minimum degree 3 (Figure \ref{C4free10-m3}). Each time we take one of the following operations, and finally we can construct a $C_4$-free planar graphs on $n$ vertices with minimum degree 3, where $10\le n\le 29$.

(1) Operation (B), as illustrated in  Figure \ref{operation(B)};

(2) The reverse operation of (B);

\begin{figure}[H]
\centering
\begin{tikzpicture}[line cap=round,line join=round,>=triangle 45,x=1.0cm,y=1.0cm]
\clip(-5.24,1.8) rectangle (6.16,4.24);
\draw [dash pattern=on 2pt off 2pt] (-5,4)-- (-5,2);
\draw (-5,4)-- (-4,4);
\draw (-4,4)-- (-3,4);
\draw (-3,4)-- (-3,2);
\draw (-5,2)-- (-4,2);
\draw (-4,2)-- (-3,2);
\draw (-3,4)-- (-2,4);
\draw (-2,4)-- (-1,4);
\draw [dash pattern=on 2pt off 2pt] (-1,4)-- (-1,2);
\draw (-1,2)-- (-2,2);
\draw (-2,2)-- (-3,2);
\draw (2,4)-- (3,4);
\draw [dash pattern=on 2pt off 2pt] (2,4)-- (2,2);
\draw (3,4)-- (4,4);
\draw (4,4)-- (4,2);
\draw (2,2)-- (3,2);
\draw (3,2)-- (4,2);
\draw (4,2)-- (5,2);
\draw (4,4)-- (5,4);
\draw (5,4)-- (6,4);
\draw [dash pattern=on 2pt off 2pt] (6,4)-- (6,2);
\draw (5,2)-- (6,2);
\draw (4,3.38)-- (2,2);
\draw (4,2.6)-- (6,4);
\draw (0.2,3.44) node[anchor=north west] {$\Longrightarrow$};
\begin{scriptsize}
\fill [color=black] (-5,4) circle (1.5pt);
\fill [color=black] (-5,2) circle (1.5pt);
\fill [color=black] (-4,4) circle (1.5pt);
\fill [color=black] (-3,4) circle (1.5pt);
\fill [color=black] (-3,2) circle (1.5pt);
\fill [color=black] (-4,2) circle (1.5pt);
\fill [color=black] (-2,4) circle (1.5pt);
\fill [color=black] (-1,4) circle (1.5pt);
\fill [color=black] (-1,2) circle (1.5pt);
\fill [color=black] (-2,2) circle (1.5pt);
\fill [color=black] (2,4) circle (1.5pt);
\fill [color=black] (3,4) circle (1.5pt);
\fill [color=black] (2,2) circle (1.5pt);
\fill [color=black] (4,4) circle (1.5pt);
\fill [color=black] (4,2) circle (1.5pt);
\fill [color=black] (3,2) circle (1.5pt);
\fill [color=black] (5,2) circle (1.5pt);
\fill [color=black] (5,4) circle (1.5pt);
\fill [color=black] (6,4) circle (1.5pt);
\fill [color=black] (6,2) circle (1.5pt);
\fill [color=black] (4,3.38) circle (1.5pt);
\fill [color=black] (4,2.6) circle (1.5pt);
\end{scriptsize}
\end{tikzpicture}
\caption{Operation (C).\label{operation(C)}}
\end{figure}
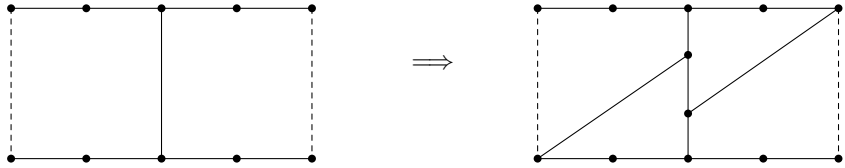

(3) Operation (C): take one edge $e$ which are the common edge of two faces of lengths at least 6, split $e$ in to three edges and add two more chordal edges between that two faces (as illustrated in Figure \ref{operation(C)}).

(iii) If $n=9$, suppose on the contrary that $\delta(9,C_4)=3$, then there is a $C_4$-free planar graph $G$ on 9 vertices with $\delta(G)=3$. By Theorem \ref{thm3}, we have $\varepsilon (G)\le 15$. Since $\delta(G)=3$, we have $\tau(G)\ne 0$, which implies that $\varepsilon(G)\le14$, and the equality holds if and only if $\tau(G)=2$ and $f_6=1$ in Theorem \ref{thm3}.  Note that the degree sequence of $G$ must be $3^84^1$,  this implies that $\tau(G)\ge 4$, since each vertex of degree 3 is adjacent of at least one edge in $\tau(G)$. This is contradicts the fact that $\tau(G)=2$. So we have $\delta(9,C_4)\le 2$. 

If $n=8$, suppose on the contrary that $\delta(8,C_4)=3$, then there is a $C_4$-free planar graph $G$ on 8 vertices with $\delta(G)=3$. By Theorem \ref{thm3}, we have $\varepsilon (G)\le 12$. Since $\delta(G)=3$, we have that $\varepsilon(G)=12$, so the degree sequence of $G$ must be $3^8$,  this implies that $\tau(G)\ge 4$, since each vertex of degree 3 is adjacent of at least one edge in $\tau(G)$. Hence by Theorem \ref{thm3} we have $\varepsilon(G)\le \frac{15}{7}(8-2)-\frac{2}{7}4<12$. This is contradicts that  $\varepsilon(G)=12$.
So we have $\delta(8,C_4)\le 2$. 

 If $n=7$, suppose on the contrary that $\delta(7,C_4)=3$, then there is a $C_4$-free planar graph $G$ on 7 vertices with $\delta(G)=3$. By Theorem \ref{thm3}, we have $\varepsilon (G)\le 10$. Since $\delta(G)=3$, we have that $\varepsilon(G)=11$, a contradiction. 
 So we have $\delta(7,C_4)\le 2$. 
 
 If $n\le 6$, suppose on the contrary that $\delta(n,C_4)=3$, then there is a $C_4$-free planar graph $G$ on n vertices with $\delta(G)=3$. On the one hand, since $\delta(G)=3$, we have $\varepsilon(G)\ge \frac{3n}{2}$; On the other hand, by Theorem \ref{thm3} we have $\varepsilon(G)\le\frac{15}{7}(n-2)$, this is impossible since $n\le 6$.
So we have $\delta(6,C_4)\le 2$ and  $\delta(5,C_4)\le 2$.

Since $C_n$ is a $C_4$-free planar graph with minimum degree 2, we therefore conclude that $\delta(n,C_4)=2$ for $5\le n\le 9$. \qed

\section{Proof of Theorem \ref{thm20111103} }

In this section we begin to consider planar Ramsey numbers of $C_4$ versus wheels. The following two Lemmas are well known.
\begin{Lemma}(Dirac) \label{Dirac}
If $\delta(G)\ge\frac{n}{2}$, then $G$ is Hamiltonian.
\end{Lemma}

\begin{Lemma}\label{lem7}(Chv\'atal-Erd\"os)
If $\alpha(G)\le k(G)$, then $G$ is Hamiltonian.
\end{Lemma}

In \cite{MR1478241}, Brandt  proved  that 
\begin{Lemma}\label{lem14}
Every non-bipartite graph of order $n$ with more than $(n-1)^2/4+1$ edges contains cycles of every length between 3 and the length of a longest cycle.
\end{Lemma}

\begin{Lemma}\label{lem15}
Let $G$ be a $C_4$-free planar graph, then its independence number $\alpha(G^c)\le 3$.
\end{Lemma}

{\bf Proof.} If $\alpha(G^c)\ge 4$, then $G$ contains a $K_4$, and hence contains a $C_4$, which contradicts the initial hypothesis. \qed

\begin{Lemma}\label{lem16}
Let $G$ be a $C_4$-free planar graph with order $n\ge 6$ and $k(G^c)\le 2$, then
there exists two vertices $x,y$ which separates some vertex $z$ from the rest in $G^c$, and further more, 
$G-\{x,y,z\}$ contains no path of length 2 in $G$.
\end{Lemma}

{\bf Proof.}  (a) Since $k(G^c)\le 2$,  there exists two vertices $x,y$ which separates $U_1$ from the rest $U_2$ in $G^c$. Then each vertex of $U_1$ is adjacent to every vertex of $U_2$ in $G$. If both $|U_1|\ge 2$ and $|U_2|\ge 2$, then $G$ will contain a $C_4$, a contradiction.

(b) Note that $z$ is adjacent to each vertex of $V(G)-\{x,y,z\}$ in $G$.  If $G-\{x,y,z\}$ contains a path of length 2 in $G$, then there will be a $C_4$ in $G$, a contradiction.
\qed

\begin{Lemma}\label{lem17}(I. Gorgol and A. Rucinski \cite{PR_Cycle})
$PR(C_4,C_3)= PR(C_3,C_4)=PR(C_4,C_5)=7$, and $PR(C_4, C_n)=n+1$ for $n\ge 6$.
\end{Lemma}

\begin{Lemma}\label{lem20111103}
Let $G$ be a $C_4$-free planar graph on $n \ge 7$ vertices, then $G^c$ contains cycles of lengths from 3 to $n-1$.
\end{Lemma}

{\bf Proof.} Let $G$ be a $C_4$-free planar graph on $n \ge 7$ vertices, by Lemma \ref{lem17}, $G^c$ contains a $C_{n-1}$. Furthermore, $G^{c}$ is not a bipartite  graph, otherwise, there will be at least one partite set with cardinality at least 4 since  $n\ge 7$, which will induce a complete graph in $G$, and hence $G$ will contain a 4-cycle, a contradiction.
By Theorem \ref{thm3}, the number of edges of $G$ is at most $\frac{15}{7}(n-2)$. So the number of edges of $G^c$ is at least ${n\choose 2}-\frac{15}{7}(n-2)> \frac{(n-1)^2}{4}+1$, for $n\ge 7$. By Lemma \ref{lem14}, $G^c$ contains cycles of lengths from 3 to $n-1$.  \qed

Bielak and Gorgol proved that $PR(C_4,K_4)=10$, since $K_4=W_3$, so $PR(C_4,W_3)=10$.

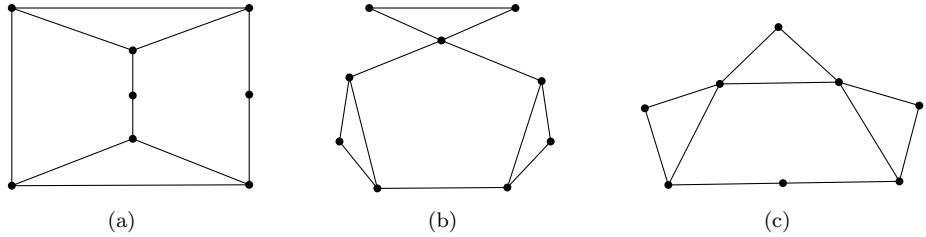
\begin{figure}[H]
\centering 
\subfigure[]{ \begin{tikzpicture}[line cap=round,line join=round,>=triangle 45,x=0.6cm,y=0.6cm]
\clip(-2.22,0.36) rectangle (3.44,4.7);
\draw (3.24,0.58)-- (3.24,2.58);
\draw (3.24,0.58)-- (-2.02,0.56);
\draw (3.24,0.58)-- (0.66,1.6);
\draw (3.24,4.5)-- (3.24,2.58);
\draw (3.24,4.5)-- (0.66,3.56);
\draw (3.24,4.5)-- (-2.02,4.5);
\draw (-2.02,4.5)-- (0.66,3.56);
\draw (-2.02,4.5)-- (-2.02,0.56);
\draw (0.66,3.56)-- (0.66,2.56);
\draw (-2.02,0.56)-- (0.66,1.6);
\draw (0.66,1.6)-- (0.66,2.56);
\begin{scriptsize}
\fill [color=black] (3.24,0.58) circle (1.5pt);
\fill [color=black] (3.24,4.5) circle (1.5pt);
\fill [color=black] (3.24,2.58) circle (1.5pt);
\fill [color=black] (-2.02,4.5) circle (1.5pt);
\fill [color=black] (0.66,3.56) circle (1.5pt);
\fill [color=black] (-2.02,0.56) circle (1.5pt);
\fill [color=black] (0.66,1.6) circle (1.5pt);
\fill [color=black] (0.66,2.56) circle (1.5pt);
\end{scriptsize}
\end{tikzpicture}}\hskip 0.4in
\subfigure[]{\begin{tikzpicture}[line cap=round,line join=round,>=triangle 45,x=0.6cm,y=0.6cm]
\clip(-0.18,0.36) rectangle (4.8,4.7);
\draw (2.18,3.78)-- (0.14,2.96);
\draw (2.18,3.78)-- (4.4,2.88);
\draw (2.18,3.78)-- (3.82,4.5);
\draw (2.18,3.78)-- (0.58,4.5);
\draw (0.76,0.5)-- (0.14,2.96);
\draw (0.76,0.5)-- (-0.08,1.54);
\draw (0.76,0.5)-- (3.64,0.52);
\draw (0.14,2.96)-- (-0.08,1.54);
\draw (3.64,0.52)-- (4.6,1.54);
\draw (4.6,1.54)-- (4.4,2.88);
\draw (3.64,0.52)-- (4.4,2.88);
\draw (3.82,4.5)-- (0.58,4.5);
\begin{scriptsize}
\fill [color=black] (2.18,3.78) circle (1.5pt);
\fill [color=black] (0.76,0.5) circle (1.5pt);
\fill [color=black] (0.14,2.96) circle (1.5pt);
\fill [color=black] (-0.08,1.54) circle (1.5pt);
\fill [color=black] (3.64,0.52) circle (1.5pt);
\fill [color=black] (4.6,1.54) circle (1.5pt);
\fill [color=black] (4.4,2.88) circle (1.5pt);
\fill [color=black] (3.82,4.5) circle (1.5pt);
\fill [color=black] (0.58,4.5) circle (1.5pt);
\end{scriptsize}
\end{tikzpicture}}\hskip 0.4in
\subfigure[]{\begin{tikzpicture}[line cap=round,line join=round,>=triangle 45,x=0.6cm,y=0.6cm]
\clip(-3.44,0.54) rectangle (3.02,4.38);
\draw (2.4,0.84)-- (2.84,2.52);
\draw (2.4,0.84)-- (1.06,3.04);
\draw (2.4,0.84)-- (-0.18,0.8);
\draw (-0.28,4.26)-- (1.06,3.04);
\draw (-0.28,4.26)-- (-1.58,3);
\draw (2.84,2.52)-- (1.06,3.04);
\draw (1.06,3.04)-- (-1.58,3);
\draw (-1.58,3)-- (-2.72,0.76);
\draw (-1.58,3)-- (-3.24,2.46);
\draw (-0.18,0.8)-- (-2.72,0.76);
\draw (-2.72,0.76)-- (-3.24,2.46);
\begin{scriptsize}
\fill [color=black] (2.4,0.84) circle (1.5pt);
\fill [color=black] (-0.28,4.26) circle (1.5pt);
\fill [color=black] (2.84,2.52) circle (1.5pt);
\fill [color=black] (1.06,3.04) circle (1.5pt);
\fill [color=black] (-1.58,3) circle (1.5pt);
\fill [color=black] (-0.18,0.8) circle (1.5pt);
\fill [color=black] (-2.72,0.76) circle (1.5pt);
\fill [color=black] (-3.24,2.46) circle (1.5pt);
\end{scriptsize}
\end{tikzpicture}}
\caption{$C_4$-free $W_n$-free planar graphs for $n=4,5,6$ respectively. \label{plg}}
\end{figure}

In Figure \ref{plg} we illustrate three $C_4$-free planar graphs which contain no $W_4, W_5$ and $W_6$ respectively, this implies that $PR(C_3,W_4)\ge 9, \ PR(C_3,W_5)\ge 10, \ PR(C_3,W_6)\ge 9$. As a matter of fact, by using a program ``Planram" due to Andrzej Dudek \cite{planram}, we can easily check the following planar ramsey numbers:  
\begin{Lemma}\label{smallram}
$PR(C_4,W_4)=9, \ PR(C_4,W_5)=10, \ PR(C_4,W_6)=9$.
\end{Lemma}

\begin{Lemma}\label{lem20111104-1}
Let $G$ be a $C_4$-free planar graph on 11 vertices, then $G^c$ contains a $W_7$.
\end{Lemma}

{\bf Proof.} By Corollary \ref{coro20101010}, we know that $\delta(G)\le 3$.   Let $v$ be a vertex in $G$ such that $d_G(v)=\delta(G)$,  let $H$ be the subgraph induced by the vertex set $V-N_G[v]$ in $G$, then $H$ is a $C_4$-free planar graph on $10-d_G(v)$ vertices. It suffices to show that $H^c$ contains a $C_7$. 

{\bf Case 1.} $\delta(G)\le 2$.

In this case we have $|V(H)|\ge 8$. Let $U\subseteq V(H)$ such that $|U|=8$, then $G[U]$ is a $C_4$-free planar graph on 8 vertices, by Lemma \ref{lem17}, $G^c[U]$ contains a $C_7$, and hence $H^c$ contains a $C_7$ too.

{\bf Case 2.} $\delta(G)=3$.

In this case we have $|V(H)|=7$.  It suffices to show that $H^c$ is Hamiltonian.

Let $t$ be the number of vertices which have degrees 3 in $G$. Since $G$ is $C_4$-free, for every vertex $u$ which has  odd degree, there must be at least one edge in $\Gamma(G)$  which is incident with $u$. This implies that $\tau(G)\ge\frac{t}{2}$. By Theorem \ref{thm3}, we have $\frac{1}{2}(3t+4(11-t))=\varepsilon(G)\le \frac{15}{7}(11-2)-\frac{2}{7}\cdot \frac{t}{2}$, this implies that $t\ge 8$.

If $t=8$, assume that  $\Delta(G)\ge 5$, then $\varepsilon(G)\ge 19$, but by 
Theorem \ref{thm3}, $\varepsilon(G)\le \frac{15}{7}(11-2)-\frac{2}{7}
\tau(G)<19$, a contradiction. So we assume that $\Delta(G)=4$, this means that the degree sequence of $G$ is $3^84^3$, hence $\varepsilon(G)=18$. On the other hand, since $\tau(G)\ge 4$, it is obvious by Theorem \ref{thm3} that $\varepsilon(G)=\frac{15}{7}(11-2)-\frac{2}{7}\tau(G)-\frac{3}{7}f_6-\frac{6}{7} f_7-\cdots -\frac{3(r-5)}{7} f_r\ne 18$ (where $r$ is the maximum length of face in $G$), a contradiction.

If $t=9$, then $\varepsilon(G)\ge 18$ and $\tau(G)\ge 5$. On the other hand, by 
Theorem \ref{thm3}, $\varepsilon(G)\le \frac{15}{7}(11-2)-\frac{2}{7}
\tau(G)<18$, a contradiction.

So the only possible case is that $t=10$. If $\Delta(G)\ge 5$, then $\varepsilon(G)\ge 18$, On the other hand, by 
Theorem \ref{thm3}, $\varepsilon(G)\le \frac{15}{7}(11-2)-\frac{2}{7}
\tau(G)<18$, a contradiction.

So we assume that $\Delta(G)=4$, and thus the degree sequence of $G$ is $3^{10}4^1$. We choose $v$ such that $d_G(v)=3$ and the only vertex of degree 4 belongs to $N_G(v)$. 

If $k(H^c)\ge 3$, then by Lemmas \ref{lem7} and \ref{lem15}, $H^c$ is Hamiltonian. So in the following, we may assume that $k(H^c)\le 2$. By Lemma \ref{lem16}, there exists two vertices $x,y$ which separates $z$ from the rest in $H^c$, which implies that $d_{H^c}(z)\le 2$, and thus $4\ge d_G(z)\ge d_H(z)\ge 4$, so $d_G(z)=4$. By the choice of $v$, we know that $z\in N_G(v)$, a contradiction.   \qed

Note that if $G$ is a planar graph of order $N$ with $\delta=\delta(G)$, then $G^c$ can not contain a $W_{N-\delta}$,  by Theorem \ref{thm20101011}, we get the following lower bounds of planar Ramsey numbers:

\begin{corollary}\label{col20120414}
\[
PR(C_4,W_n)\ge\left\{
\begin{array}{ll}
n+4, & {\rm if}\ n\in\{k|7\le k\le 25\}\cup\{27,28,29,30,31,33,34,36,37,39\};\\
\\
n+5, & {\rm if}\ n\in \{26,32,35,38\}\cup\{k|k\ge 40\}.
\end{array}
\right.
\]

\end{corollary}

\begin{Lemma}\label{lem18}
Let $G$ be a $C_4$-free planar graph on $N\ (N\ge 12)$ vertices and let $n=N-\delta(N,C_4)-1$, then $G^c$ contains  $W_n$ and $W_{n-1}$.
\end{Lemma}

{\bf Proof.} By Corollary \ref{coro20101010}, we know that $\delta(G)\le\delta(N,C_4)\le 4$. Since $n=N-\delta(N,C_4)-1$ and $N\ge 12$, we get that $n\ge 7$.  Let $v$ be a vertex in $G$ such that $d_G(v)=\delta(G)$,  let $H$ be the subgraph induced by the vertex set $V-N_G[v]$ in $G$, then $H$ is a $C_4$-free planar graph on $N-d_G(v)-1\ge n$ vertices. 

{\bf Case 1.} $\delta(G)\le \delta(N,C_4)-1$.

In this case we have  $|V(H)|\ge n+1$,  by Lemma \ref{lem20111103}, $H^c$ contains cycles of lengths from 3 to $|V(H)|-1\ge n$, let $C_{n-1}$ and $C_n$ be the cycles of lengths $n-1$ and $n$ respectively, hence $v+C_{n-1}$  and $v+C_{n}$ are  $W_{n-1}$ and $W_n$  in $G^c$ respectively.

{\bf Case 2.} $\delta(G) =\delta(N,C_4)$. 

In this case we have $|V(H)|=n$.

By Lemma \ref{lem17}, $H^c$ contains a $C_{n-1}$, and hence $v+C_{n-1}$ is a $W_{n-1}$ in $H^c$. Next, we shall show that $H^c$ contains a $W_n$.

If $k(H^c)\ge 3$, then by Lemmas \ref{lem7} and \ref{lem15}, $H^c$ is Hamiltonian. Let $C$ be a Hamiltonian cycle in $H^c$, then $v+C$ is a $W_n$ in $G^c$. So in the following, we may assume that $k(H^c)\le 2$. By Lemma \ref{lem16}, there exists two vertices $x,y$ which separates $z$ from the rest. Let $U=V(H)-\{x,y,z\}$.

Note that $z$ is adjacent to each vertex of $U$ in $G$.

If $\delta(G)=\delta(N,C_4)=4$, the number of edges of $G$
is at least $\frac{1}{2}((n-3)+4(n+4))\ge \frac{15}{7}(n+3)+1=\frac{15}{7}(|V(G)|-2)+1$, which contradicts Theorem \ref{thm3}. 

Since $N\ge 12$,  we assume that $\delta(G)=\delta(N,C_4)=3$ by  Theorem \ref{thm20101011}. In this case $|U|=N-7$.

Since $G$ is $C_4$-free and $z$ is adjacent to each vertex of $U$ in $G$, each vertex of $N_G(v)\cup\{x,y\}$ can be adjacent to at most one vertex in $V(H)-\{x,y,z\}$ in $G$. So the edges between $N_G[v]\cup\{x,y,z\}$ and $U$ is at most 
$|U|+5$; On the other hand, since $\delta(G)=3$ and their is no path of length 2 in the subgraph of $G$ induced by $U$ by Lemma \ref{lem16}, the number of edges between $N_G[v]\cup\{x,y,z\}$ and $U$ is at least $3|U|-2[\frac{|U|}{2}]$ (where $[x]$ denotes the maximum integer which is at most $x$). So we have that $3|U|-2[\frac{|U|}{2}]\le |U|+5$, which is impossible since $N\ge 12$.\qed

Combining Theorem \ref{thm20101011}, corollary \ref{col20120414} and Lemmas \ref{lem20111104-1},\ref{smallram} and \ref{lem18}, we finally prove Theorem \ref{thm20111103}.


\begin{thebibliography}{10}
\bibitem{MR1478241}
S. Brandt.
\newblock A sufficient condition for all short cycles.
\newblock {\em Discrete Appl. Math.}, 79(1-3):63--66, 1997.
\newblock 4th Twente Workshop on Graphs and Combinatorial Optimization
  (Enschede, 1995).
\bibitem{dudek}
A. Dudek and A. Rucinski, Planar Ramsey numbers for small graphs, 36th South- eastern International Conference on Combinatorics, Graph Theory, and Computing. Congr. Numer. 176 (2005), 201220.
\bibitem{PR_Cycle} I. Gorgol and A. Rucinski,
\newblock Planar ramsey numbers for cycles,
\newblock  Discrete Mathematics,  308(2008), 4389-4395.

\bibitem{MR1244935}
R. Steinberg and C.~A. Tovey,
\newblock Planar {R}amsey numbers,
\newblock  Journal of Combinatorial Theory (B), 59(1993), 288--296.

\bibitem{walker} K. Walker, The analog of Ramsey numbers for planar graphs, Bulletin of London Mathematics, 1(1967), 187-190.

\bibitem{zhouC3} Zhou Guofei, Chen Yaojun, Miao Zhengke and S. Pirzada, 
A note on planar Ramsey numbers for a triangle versus wheels,
 Discrete Mathematics and Theoretical Computer Science
 14(2)  (2012), 255-260. 
 
\bibitem{chen} Zhou Guofei and Chen Yaojun, The Maximum Size of $C_4$-free Planar Graphs, accepted to be published in Ars Combinatoria.

\bibitem{plantri} http://cs.anu.edu.au/~bdm/plantri

\bibitem{planram}
A. Dudek, http://homepages.wmich.edu/~zyb1431/papers/planram.c

  
\end{thebibliography}
\end{document}